\documentclass[12pt,reqno]{amsproc}
\usepackage{amssymb}

\usepackage{float}
\usepackage{amsmath}
\usepackage{graphicx}
\usepackage{amscd}

\theoremstyle{plain}
\newtheorem{X}{X}[section]
\newtheorem{theorem}[X]{Theorem}
\newtheorem{proposition}[X]{Proposition}
\newtheorem{lemma}[X]{Lemma}
\newtheorem{corollary}[X]{Corollary}
\newtheorem{conjecture}[X]{Conjecture}

\theoremstyle{definition}

\theoremstyle{remark}

\newtheorem{remark}[X]{Remark}
\newtheorem*{Remark}{Remark}

\numberwithin{equation}{section}
\numberwithin{equation}{section}

\usepackage{verbatim}
\usepackage[size=2em,heads=vee]{diagrams}
\diagramstyle[labelstyle=\scriptstyle]

\makeindex

\newenvironment{mylist}
{\begin{list}
{--}
{\setlength{\leftmargin}{.5in}\setlength{\rightmargin}{.5in}}}
{\end{list}}

\renewenvironment{itemize}
{\begin{mylist}}
{\end{mylist}}

\numberwithin{equation}{section}

\let\hat=\widehat  
\let\tilde=\widetilde  
\def\1{{1\mkern-7mu1}}  

\newcommand\Aut{\operatorname{Aut}}

\newcommand\codim{\operatorname{codim}}

\newcommand\Hom{\operatorname{Hom}}

\newcommand\Ker{\operatorname{Ker}}

\newcommand\Lie{\operatorname{Lie}}

\let\plim=\varprojlim

\newcommand\rank{\operatorname{rank}}

\newcommand\Res{\operatorname{Res}}

\def\tint{\mathop{\textstyle \int}}%

\begin{document}
\title{Kazhdan's Theorem on Arithmetic Varieties}
\dedicatory{}
\author{J.S. Milne}

\begin{abstract}
Define an arithmetic variety to be the quotient of a bounded symmetric
domain by an arithmetic group. An arithmetic variety is algebraic, and the
theorem in question states that when one applies an automorphism of the
field of complex numbers to the coefficients of an arithmetic variety the
resulting variety is again arithmetic. This article simplifies Kazhdan's proof. In particular, it avoids recourse to the
classification theorems.

It was originally completed on March 28, 1984, and distributed in
handwritten form.
\end{abstract}

\maketitle
\tableofcontents

\setcounter{section}{0}

Let $X$ be a nonsingular algebraic variety over $\mathbb{C}$ whose universal
covering manifold $\tilde {X}$ is a symmetric Hermitian domain (i.e., a
symmetric Hermitian space without Euclidean or compact factors). Then the
group $\Aut (\tilde {X})$ of automorphisms of $\tilde {X}$ (as a complex
manifold) has only finitely many connected components, and its identity
component $\Aut(\tilde {X} )^{+}$ is a real semisimple Lie group with no
compact factors. The variety $X$ will be said to be \emph{arithmetic\/} if
the group $\Gamma$ of covering transformations of $\tilde {X}$ over $X$ is
torsion-free and is an arithmetic subgroup of $\Aut(\tilde {X})$ in the
sense that there is a linear algebraic group $G$ over $\mathbb{Q}$ and a
surjective homomorphism $f\colon G(\mathbb{R})^{+}\to\Aut(\tilde {X})^{+}$ with
compact kernel carrying $G (\mathbb{Z})\cap G(\mathbb{R})^{+}$ into a group
commensurable with $\Gamma$.

The theorem in question is the following:

\begin{theorem}
If $X$ is an arithmetic variety, then, for all automorphisms $\sigma$ of $%
\mathbb{C}$, $\sigma X$ is also an arithmetic variety.
\end{theorem}

The result was first stated, with some indications of a proof, by Kazhdan in
his talk at the Nice Congress in 1970 (Kazhdan 1971), but, in fact, it is
only recently\footnote{Added 22.06.01. Recall that the article was written in 1984.} that the proof has been completed.\footnote{%
Added 22.06.01. For a different proof, see
\begin{quotation}
Nori, M. V.; Raghunathan, M. S., On conjugation of locally symmetric
arithmetic varieties. Proceedings of the Indo-French Conference on Geometry
(Bombay, 1989), 111--122, Hindustan Book Agency, Delhi, 1993.
\end{quotation}}

In the case that $X$ is compact, a detailed proof of (0.1) was given by
Kazhdan in his talk at the Budapest conference in 1971 (Kazhdan 1975). There
is now an alternative approach to this case of the theorem: roughly
speaking, the first part of Kazhdan's paper, where the nondegeneracy of the
Bergmann metric on the universal covering manifold of $\sigma X$ is
established, can be replaced by an appeal to Yau's theorem on the existence
of K\"{a}hler-Einstein metrics; the second, group-theoretic part, of the
paper can be replaced by an appeal to the general results of Margulis (see
(2.7) below).

Let $X^{\prime }\rightarrow X$ be a finite \'{e}tale morphism. If $X$ is
arithmetic, then clearly so also is $X^{\prime }$, and the converse
assertion is true provided the fundamental group of $X$ is torsion-free.
From this it follows that it suffices to prove (0.1) under the hypothesis:

\begin{quote}
(0.2) There exists an almost $\mathbb{Q}$-simple, simply connected algebraic
group $G$ over $\mathbb{Q}$, a surjective homomorphism $f\colon G(\mathbb{R}%
)\rightarrow \Aut(\tilde{X})^{+}$ with compact kernel, and a congruence
subgroup $\Gamma \subset G(\mathbb{Q})$ such that $X$ is equal to $f(\Gamma
)\backslash \tilde{X}$ with its unique algebraic structure.
\end{quote}

\noindent The group $G$ will then be of the form $G=\Res_{F/\mathbb{Q}%
}G^{\prime }$ with $G^{\prime }$ an absolutely almost-simple group over a
totally real number field $F$. The type of $G_{\mathbb{R}}$ ($A,B,C,D^{%
\mathbb{R}},D^{\mathbb{H}},$ mixed type $(D^{\mathbb{R}},D^{\mathbb{H}})$, $%
E_{6}$, $E_{7}$; cf. Deligne 1979) will be called the \emph{type\/} of $X$.

If $X$ is as in (0.2) and of type $A,B,C,D^{\mathbb{R}}$, or $D^{\mathbb{H}}$%
, then it is a moduli variety for a class of abelian varieties with a family
of Hodge cycles and a level structure---see for example (Milne and Shih
1982, \S 2). (In the case $D^{\mathbb{H}}$ one should assume that the image
of $\Gamma $ in $G^{\text{ad}}(\mathbb{Q})$ is a congruence subgroup there.)
For any automorphism $\sigma $ of $\mathbb{C}$, the conjugate variety $%
\sigma X$ will be a moduli variety for the conjugate class of abelian
varieties with extra structure, and so is again an arithmetic variety. Thus
(0.1) is true for arithmetic varieties of these types. (Note that this
argument makes use of Deligne's theorem (Deligne 1982) that a Hodge cycle on
an abelian variety is absolutely Hodge.)

Kazhdan (1983) proves the following result:

\begin{quote}
let $X$ be as in (0.2), noncompact, and of type $E_6$, $E_7$, or mixed type $%
(D^{\mathbb{R}},D^{\mathbb{H}})$; then, if one assumes a certain result in
the theory of group representations (Conjecture 3.17 below), (0.1) is true
for $X$.
\end{quote}

\noindent A weak form of (3.17), sufficient for the purpose of proving
(0.1), has recently been established by L. Clozel.\footnote{Added 22.06.01. Clozel, Laurent,
On limit multiplicities of discrete series representations in spaces of automorphic forms.
Invent. Math. 83 (1986), no. 2, 265--284.

See also:

Rohlfs, J., and Speh, B., On the limit multiplicities of representations with cohomology in
the cuspidal spectrum. Duke Math. J. 55 (1987), no. 1, 199--211.} Thus this last result of
Kazhdan, together with the proofs for compact varieties and moduli
varieties, suffices to prove (0.1).

For his proof in these last cases, Kazhdan shows, case-by-case, that $X$
contains an arithmetic subvariety $X^{\prime }$ associated with a subgroup $%
G^{\prime }$ of $G$ such that 
\begin{equation*}
\dim X>\dim X^{\prime }>\dim (X^{\ast }\smallsetminus X)
\end{equation*}%
where $X^{\ast }$ is the canonical (Baily-Borel) compactification of $X$.
Because (0.1) is known for curves (they are all of type $A_{1}$), in proving
(0.1) for $X$ he can assume by induction on the dimension that it holds for $%
X^{\prime }$.

In these notes, I modify Kazhdan's arguments to give a uniform proof that
(0.1) holds for all $X$ (satisfying (0.2)) such that

\begin{enumerate}
\item[(0.3a)] \ $\codim(X^{\ast }\smallsetminus X)\geq 3$ (i.e., the codimension of
the boundary in $X^{\ast}$ is at least $3$);

\item[(0.3b)] \ $G^{\prime }$ contains a maximal torus $T^{\prime }$ that
splits over an imaginary quadratic extension of $F$ and is such that $%
T^{\prime }(\mathbb{R})$ is compact.
\end{enumerate}

\noindent The proof does not use the classification and avoids treating the
compact varieties separately, but it does assume (0.1) for the case of
arithmetic varieties of type $A_{1}$ (for these varieties it is easily
proved using moduli varieties). It also, of course, makes use of Clozel's
result.

Let $X$ satisfy (0.2). A construction of Piatetsky-Shapiro and Borovoi
allows one to embed $X$ into an arithmetic variety $X^{\prime}$ satisfying
also (0.3). I do not know whether it is possible to prove directly that
(0.1) for $X^{\prime}$ implies (0.1) for $X$, but in any case one has the
implications:

\begin{quote}
(0.1) for all arithmetic varieties satisfying (0.2) and (0.3) $\implies$
Langlands's conjecture (Langlands 1979, p232--33) for all Shimura varieties $%
\implies$ (0.1) for all arithmetic varieties.
\end{quote}

\noindent See Milne 1983 or Borovoi ...\footnote{Added 11.07.01: 
Borovoi, M. V., Langlands' conjecture concerning conjugation of connected Shimura varieties. Selected translations. Selecta Math. Soviet. 3 (1983/84), no. 1, 3--39.

Borovoi, M. V., On the groups of points of a semisimple group over a totally real field, Problems in Group
Theory and Homological Algebra, Yaroslavl 1987, pp.~142--149}

The first section of these notes contains generalities on arithmetic
varieties. Also a group $Q$ is defined that will play the role of $G(\mathbb{%
Q})/Z(\mathbb{Q})$ for the conjugate variety $\sigma X$. Criteria are given
in \S 2 for a variety to be arithmetic; in conjunction with Yau's theorem,
they suffice to show that $\sigma X$ is arithmetic in the case that $X$ is
compact. In \S 3 Clozel's result is used to show that the Bergmann volume
form on a certain covering manifold of $\sigma X$ is not identically zero.
In the next section, the subbundles of the tangent bundle on $\tilde{X}$ are
studied and, finally, in \S 5 the main theorem is proved.

\section{Generalities on arithmetic varieties; definition of $\tilde{X}^{%
\protect\sigma }$ and $Q$}

Let $G$ be a simply connected semisimple algebraic group over $\mathbb{Q}$,
and let $\tilde{X}$ be a symmetric Hermitian domain on which $G(\mathbb{R})$
acts in such a way that
\[
\tilde{X}\approx G(\mathbb{R})/\text{(maximal compact subgroup).}
\]
We always assume that $G$ has no $\mathbb{Q}{}$-factor that is
anisotropic over $\mathbb{R}{}$. For any compact open subgroup $K$ of $G(%
\mathbb{A}_{f})$, $K\cap G(\mathbb{Q})$ is a congruence subgroup of $G(%
\mathbb{Q})$, and we let $\Gamma =\Gamma (K)$ be its image in $G^{\text{ad}}(%
\mathbb{Q})$. Clearly $\Gamma (K)=\Gamma (K\cdot Z(\mathbb{Q}))$, where $%
Z=Z(G)$, and so we can always assume that $K\supset Z(\mathbb{Q}{})$. Let $%
X_{K}$, or $X_{\Gamma }$, denote $\Gamma (K)\backslash \tilde{X}$ regarded
as an algebraic variety. The strong approximation theorem shows that $G(%
\mathbb{Q}{})$ is dense in $G(\mathbb{A}_{f})$, and it follows that 
\begin{equation*}
X_{K}(\mathbb{C}{})\overset{\text{df}}{=}\Gamma \backslash \tilde{X}\cong G(%
\mathbb{Q})\backslash \tilde{X}\times G(\mathbb{A}_{f})/K.
\end{equation*}

\noindent The actions on the right hand term are 
\begin{equation*}
q(x,a)k=(qx,qak),\quad q\in G(\mathbb{Q}{}),\quad x\in \tilde{X},\quad a\in
G(\mathbb{A}_{f}),\quad k\in K.
\end{equation*}

Let $\mathcal{S}$ denote the set of compact open subgroups $K$ of $G(\mathbb{%
A}_{f})$ containing $Z(\mathbb{Q})$ and such that $\Gamma (K)$ is
torsion-free. The groups $\Gamma (K)$, $K\in \mathcal{S}$, have the
following properties:

\begin{itemize}
\item each $\Gamma (K)$ is an arithmetic subgroup;

\item $\bigcap_{K\in \mathcal{S}}\Gamma (K)=\{1\}$;

\item if $K_{1},\ldots ,K_{m}$ are in $\mathcal{S}$, then $K=_{\text{df}%
}\bigcap K_{i}$ is in $\mathcal{S}$ and $\Gamma (K)=\bigcap \Gamma (K_{i})$;

\item if $q\in G^{\text{ad}}(\mathbb{Q})$ and $K\in \mathcal{S}$, $q\Gamma
(K)q^{-1}=\Gamma (qKq^{-1})$ with $qKq^{-1}\in \mathcal{S}$.
\end{itemize}

\noindent We shall often identify the elements of $\mathcal{S}$ with their
images in $G(\mathbb{A}_{f})/Z(\mathbb{Q})$.

Consider the projective system $(X_{K})_{K\in \mathcal{S}}$ of algebraic
varieties. There is a left action of $G(\mathbb{A}_{f})$ on this system: 
\begin{equation*}
g\colon X_{K}\rightarrow X_{gKg^{-1}},\quad \lbrack x,a]\mapsto \lbrack
x,ag^{-1}],\quad x\in \tilde{X},\quad a,g\in G(\mathbb{A}_{f}).
\end{equation*}

\noindent Let $\hat{X}=\plim X_{K}$ --- this is an irreducible scheme over $%
\mathbb{C}$ (not of finite-type!).

\begin{lemma}
There are natural bijections (of sets) 
\begin{equation*}
\hat{X}(\mathbb{C}{})\cong \plim X_{K}(\mathbb{C}{})\cong G(\mathbb{Q}%
{})\backslash \tilde{X}\times G(\mathbb{A}_{f})\cong (G(\mathbb{Q}{})/Z(%
\mathbb{Q}{}))\backslash \tilde{X}\times (G(\mathbb{A}_{f})/Z(\mathbb{Q}{})).
\end{equation*}
\end{lemma}

\begin{proof}
Only the middle bijection requires proof. Recall the following result
(Bourbaki 1989, III 7.2): consider a projective system $(G_{\alpha })$ of
topologicaly groups acting continuously and compatibly on a topological
space $S$; then the canonical map $S/\plim G_{\alpha }\rightarrow \plim%
(S/G_{\alpha })$ is bijective if

\begin{enumerate}
\item the isotropy group in $G_{\alpha }$ of each $s\in S$ is compact, and

\item the orbit $G_{\alpha }s$ of each $s\in S$ is compact.
\end{enumerate}

\noindent Apply this with $S=G(\mathbb{Q})\backslash \tilde{X}\times G(%
\mathbb{A}_{f})$ and $(G_{\alpha })=(K)$. Then (a) holds because each
isotropy group is $Z(\mathbb{Q})$, and (b) holds because $K$ is compact. As $%
\plim K=\cap K=1$, we have that 
\begin{equation*}
G(\mathbb{Q}{})\backslash \tilde{X}\times G(\mathbb{A}_{f})\rightarrow \plim %
G(\mathbb{Q})\backslash \tilde{X}\times G(\mathbb{A}_{f})/K
\end{equation*}

\noindent is bijective, as required.
\end{proof}

The action of $G(\mathbb{A}_{f})$ on the system $(X_{K})$ defines an action
of $G(\mathbb{A}_{f})$ on $\hat{X}$: for $g\in G(\mathbb{A}_{f})$, $g\colon 
\hat{X}\rightarrow \hat{X}$ is $[x,a]\mapsto \lbrack x,ag^{-1}]$. The
knowledge of $\hat{X}$ with this action is equivalent to the knowledge of
the projective system $(X_{K})$ together with the action of $G(\mathbb{A}%
_{f})$ on it; in particular, $X_{K}=K\backslash \hat{X}$. We shall sometimes
use $\hat{X}$ to denote the projective system $(X_{K})$ rather than its
limit. Note that the action of $G(\mathbb{A}_{f})/Z(\mathbb{Q})$ on $\hat{X}$
is effective.

There is an action of $G^{\text{ad}}(\mathbb{Q}{})^{+}$ on $\hat{X}$: 
\begin{equation*}
\alpha \lbrack x,a]=[\alpha x,\alpha (a)],\quad \alpha \in G^{\text{ad}%
}(X)^{+},\quad x\in \tilde{X},\quad a\in G(\mathbb{A}_{f}).
\end{equation*}

\noindent The actions of $G(\mathbb{Q}{})$ on $\hat{X}$ defined by the maps $%
G(\mathbb{Q})\rightarrow G^{\text{ad}}(\mathbb{Q}{})$ and $G(\mathbb{Q}%
{})\rightarrow G(\mathbb{A}_{f})$ are equal: 
\begin{equation*}
\lbrack qx,qaq^{-1}]=[x,aq^{-1}]\text{,}\quad q\in G(\mathbb{Q})\text{,}%
\quad x\in \tilde{X},\quad a\in G(\mathbb{A}_{f}).
\end{equation*}

\bigskip

Let $(X_{\alpha })$ be a projective system of smooth complex algebraic
varieties such that the transition maps $X_{\beta }\rightarrow X_{\alpha }$, 
$\beta \geq \alpha $, are \'{e}tale. Then $\hat{X}=_{\text{df}}\plim %
X_{\alpha }$ has a canonical structure as a complex manifold: a basis for
the atlas on $\hat{X}$ is formed by the pairs $(U,\varphi )$ for which there
exists an $\alpha $ such that the projection $p_{\alpha }\colon \hat{X}%
\rightarrow X_{\alpha }$ is injective on $U$ and $(p_{\alpha }(U),\varphi
\circ p_{\alpha }^{-1})$ is an open chart for $X_{\alpha }$. We write $\hat{X%
}^{\text{an}}$ for $\hat{X}$ with this structure. Note that the topology on $%
\hat{X}^{\text{an}}$ is, in general, strictly finer than the projective
limit of the topologies on the $X_{\alpha }^{\text{an}}$. (However, the
Zariski topology on $\hat{X}$ \emph{is }the projective limit of the Zariski
topologies on the $X_{\alpha }$.)

In the situation considered above, if $G(\mathbb{A}_{f})$ is given the
discrete topology, then 
\begin{equation*}
\hat{X}^{\text{an}}=G(\mathbb{Q}{})\backslash \tilde{X}\times G(\mathbb{A}%
_{f}).
\end{equation*}

\noindent Note that the map 
\begin{equation*}
\tilde{X}\rightarrow \hat{X}^{\text{an}},\quad x\mapsto \lbrack x,1]
\end{equation*}

\noindent is injective (because $\cap \Gamma (K)=\{1\}$), and is an
isomorphism of $\tilde{X}$ (as a complex manifold) onto a connected
component of $\hat{X}^{\text{an}}$. We shall use the following notations for
the maps:

\begin{diagram}
\tilde{X} &  & \rInto^{p} &  & \hat{X}^{\text{an}} &  &  &  & \hat{X} \\ 
& \rdTo_{p_{K}} &  & \ldTo_{q_{K}^{\text{an}}} &  &  &  & \ldTo_{q_{K}} &  \\ 
&  & X_{K}^{\text{an}} &  &  &  & X_{K} &  & 
\end{diagram}

When $x\in \tilde{X}$, we often write $x_{K}$ and $\hat{x}$ for $p_{K}(x)$
and $p(x)$.

\begin{lemma}
Let $g\in G(\mathbb{A}_{f})/Z(\mathbb{Q}{})$; if $g$ stabilizes $\tilde{X}%
\subset \hat{X}$, then it belongs to $G(\mathbb{Q})/Z(\mathbb{Q})$.
\end{lemma}

\begin{proof}
Let $x\in \tilde{X}$. As $g$ stabilizes $\tilde{X}$, $[x,g]=[x^{\prime },1]$
for some $x^{\prime }\in \tilde{X}$. This means that there exists a $q\in G(%
\mathbb{Q}{})/Z(\mathbb{\mathbb{Q}{})}$ such that $(qx^{\prime },q)=(x,g)$
as elements of $\tilde{X}\times G(\mathbb{A}_{f})/Z(\mathbb{Q})$. In
particular, $q=g\in G(\mathbb{Q}{})/Z(\mathbb{Q}{})$.
\end{proof}

Now fix an automorphism $\sigma $ of $\mathbb{C}{}$. The discussion above
shows that there is a canonical structure of a complex manifold on $\sigma 
\hat{X}$. Choose a connected component $\tilde{X}^{\sigma }$ of $(\sigma 
\hat{X})^{\text{an}}$, and let $p^{\sigma }$ and $p_{K}^{\sigma }$ be the
inclusion $\tilde{X}^{\sigma }\hookrightarrow (\sigma \hat{X})^{\text{an}}$
and the composite $(\sigma q_{K})^{\text{an }}\circ p^{\sigma }$%
respectively: thus 
\begin{diagram}
\tilde{X}^{\sigma } & \rInto{p^{\sigma }} & (\sigma \hat{X})^{\text{an}} &  & \sigma 
\hat{X}   \\ 
& \rdTo_{p_{K}^{\sigma }} & \dTo_{(\sigma q_{K})^{\text{an}}} &  & \dTo_{\sigma q_{K}}
\\ 
&  & (\sigma X_{K})^{\text{an}} &  & \sigma X_{K}
\end{diagram}The group $G(\mathbb{A}_{f})$ continues to act on $\sigma \hat{X%
}$, and this action is compatible with the complex structure on $(\sigma 
\hat{X})^{\text{an}}$. Define $Q\subset G(\mathbb{A}_{f})/Z(\mathbb{Q}{})$
to be the stabilizer of $\tilde{X}^{\sigma }$. Note that $\tilde{X}^{\sigma
}\rightarrow \sigma X_{K}^{\text{an}}$ is a local isomorphism. Let $x\in \tilde{X%
}^{\sigma }$, and let $M$ be the universal covering manifold of $\sigma X_{K}^{%
\text{an}}$. For any $m\in M$, there is a unique map $M\rightarrow \sigma 
\hat{X}$ such that $m$ maps to $\hat{x}$ and the composites $M\rightarrow
\sigma \hat{X}\rightarrow \sigma X_{K^{\prime }}^{\text{an}}$ are all
analytic. Clearly, $M\rightarrow \sigma \hat{X}$ is analytic, and so its
image is contained in $\tilde{X}^{\sigma }$. This shows that $\tilde{X}%
^{\sigma }\rightarrow X^{\text{an}}$ is surjective and that $\tilde{X}%
^{\sigma }$ is a covering manifold of $X^{\text{an}}$. In particular, $%
\tilde{X}^{\sigma }$ is Zariski dense in $\sigma \hat{X}$, and so the action
of $Q$ on $\tilde{X}^{\sigma }$ is effective.

\begin{lemma}
For any $K$ in $\mathcal{S}$, the map 
\begin{equation*}
(k,x)\mapsto kx\colon K\times \tilde{X}^{\sigma }\rightarrow (\sigma \hat{X}%
)^{\text{an}}
\end{equation*}
is surjective.
\end{lemma}
\begin{proof}
Let $x\in \sigma \hat{X}^{\text{an}}$. By what we have just proved, there
exists an $\tilde{x}\in \tilde{X}^{\sigma }$ such that $\tilde{x}$ and $x$
have the same image in $\sigma X_{K}^{\text{an}}$. Let $\hat{x}=p^{\sigma }(%
\tilde{x})$. Then $(\sigma q_{K})(\hat{x})=(\sigma q_{K})(x)$, and so $%
q_{K}(\sigma ^{-1}\hat{x})=q_{K}(\sigma ^{-1}x)$ in $X_{K}^{\text{an}%
}=K\backslash \hat{X}^{\text{an}}.$ Therefore, there exists a $k\in K$ such
that $k(\sigma ^{-1}\hat{x})=\sigma ^{-1}x$, and so $k\hat{x}=x$.
\end{proof}

\begin{proposition}
For any $K$ in $\mathcal{S}$, $KQ=G(\mathbb{A}_{f})/Z(\mathbb{Q})$. In
particular, $Q$ is dense in $G(\mathbb{A}_{f})/Z(\mathbb{Q})$.
\end{proposition}

\begin{proof}
For any $g\in G(\mathbb{A}_{f})$, $g\tilde{X}^{\sigma }$ is a connected
component of $\sigma \hat{X}^{\text{an}}$, and so the lemma shows there
exists a $k$ in $K$ such that $kg\tilde{X}^{\sigma }=\tilde{X}^{\sigma }$.
By definition, $kg\in Q$, and therefore $g\in k^{-1}Q\subset KQ$.
\end{proof}

\begin{proposition}
(a) The map 
\begin{equation*}
\lbrack x,g]\mapsto gx\colon Q\backslash \tilde{X}^{\sigma }\times G(\mathbb{%
A}_{f})\rightarrow \sigma \hat{X}^{\text{an}}
\end{equation*}

\noindent is an isomorphism of complex manifolds ($G(\mathbb{A}_{f})$ with
the discrete topology).

(b) For every $K$, the map 
\begin{equation*}
\lbrack x]\mapsto p_{K}^{\sigma }(x)\colon (Q\cap K)\backslash \tilde{X}%
^{\sigma }\rightarrow \sigma X_{K}^{\text{an}}
\end{equation*}

\noindent is an isomorphism of complex manifolds.
\end{proposition}

\begin{proof}
(a) The lemma shows that $\tilde{X}^{\sigma }\times G(\mathbb{A}%
_{f})\rightarrow \sigma \hat{X}^{\text{an}}$ is surjective, and it is
obvious that the fibres of the map are the orbits of $Q$.

(b) As $\sigma X_{K}=K\backslash \sigma \hat{X}$, it follows from (a) that
there is an isomorphism 
\begin{equation*}
Q\backslash \tilde{X}^{\sigma }\times G(\mathbb{A}_{f})/K\rightarrow \sigma
X_{K}^{\text{an}}.
\end{equation*}
The preceding proposition can be used to show that 
\begin{equation*}
(Q\cap K)\backslash \tilde{X}^{\sigma }\rightarrow Q\backslash \tilde{X}%
^{\sigma }\times G(\mathbb{A}_{f})/K
\end{equation*}
is an isomorphism.
\end{proof}

We write $\Gamma ^{\sigma }(K)$ for $Q\cap K$. Thus, the proposition shows
that 
\begin{equation*}
\Gamma ^{\sigma }(K)\backslash \tilde{X}^{\sigma }\overset{\cong }{%
\rightarrow }\sigma (\Gamma (K)\backslash \tilde{X}).
\end{equation*}

\begin{proposition}
For any $K\in \mathcal{S}$, $Q$ is contained in the commensurability group
of $\Gamma ^{\sigma }(K)$.
\end{proposition}

\begin{proof}
Let $g\in Q$; we have to show that $g\Gamma ^{\sigma }(K)g^{-1}$is
commensurable with $\Gamma ^{\sigma }(K)$. Note that $g\Gamma ^{\sigma
}(K)g^{-1}=\Gamma ^{\sigma }(gKg^{-1})$. Let $K^{\prime }=K\cap gKg^{-1}$.
The diagram of finite \'{e}tale maps 
\begin{diagram}
&  & X_{K^{\prime }} &  &  \\ 
& \ldTo &  & \rdTo &  \\ 
X_{K} &  &  &  & X_{gKg^{-1}}
\end{diagram}

\noindent gives rise to a similar diagram

\begin{diagram}
&  & \sigma X_{K^{\prime }}=\Gamma ^{\sigma }(K^{\prime })\backslash \tilde{X%
}^{\sigma } &  &  \\ 
& \ldTo &  &  \rdTo &\\ 
\sigma X_{K}=\Gamma ^{\sigma }(K)\backslash \tilde{X}^{\sigma } &  &  &  & 
\sigma X_{gKg^{-1}}=\Gamma ^{\sigma }(gKg^{-1})\backslash \tilde{X}^{\sigma }
\end{diagram}
This shows that $\Gamma ^{\sigma }(K^{\prime })$ has finite index in both $%
\Gamma ^{\sigma }(K)$ and $\Gamma ^{\sigma }(gKg^{-1})=g\Gamma ^{\sigma
}(K)g^{-1}$.
\end{proof}

\begin{remark}
It follows from (1.3) that $G(\mathbb{A}_{f})$ acts transitively on the
space of connected components of $(\sigma \hat{X})^{\text{an}}$. Thus every
connected component is of the form $g\tilde{X}^{\sigma }$, $g\in G(\mathbb{A}%
_{f}),$ and has a stabilizer that is conjugate to $Q$.
\end{remark}

\begin{proposition}
Let $Z$ be a nonempty $Q$-invariant subset of $\tilde{X}^{\sigma }$. Assume
that $Z$ is a complex analytic subset of an open submanifold of $\tilde{X}%
^{\sigma }$ and that $p_{K}^{\sigma }(Z)$ is an algebraic subvariety of $%
\sigma X_{K}$ for some $K$ in $\mathcal{S}$; then $Z=\tilde{X}^{\sigma }$.
\end{proposition}

\begin{proof}
The hypotheses imply that $\hat{Z}=_{\text{df}}(p_{K}^{\sigma }(Z))_{K\in 
\mathcal{S}}$ is a pro-algebraic subvariety of $\sigma \hat{X}$ and is $Q$%
-invariant. As $Q$ is dense in $G(\mathbb{A}_{f})/Z(\mathbb{Q}{})$, $\hat{Z}$
is invariant under $G(\mathbb{A}_{f})/Z(\mathbb{Q})$, and so $\sigma ^{-1}(%
\hat{Z})$ is a $G(\mathbb{A}_{f})/Z(\mathbb{Q})$-invariant subscheme of $%
\hat{X}$. Let $\tilde{Z}=\tilde{X}\cap \sigma ^{-1}(\hat{Z})$. Then $\tilde{Z%
}$ is $G(\mathbb{Q})$-invariant and, as $G(\mathbb{Q})$ is dense in $G(%
\mathbb{R})$ and $G(\mathbb{R})$ acts transitively, this shows that $\tilde{Z%
}=\tilde{X}$.
\end{proof}

For any $K$ in $\mathcal{S}$, let $X_{K}^{\ast }$ be the canonical
(Baily-Borel) compactification of $X_{K}$. The dimension of $X_{K}^{\ast
}\smallsetminus X_{K}$ is independent of $K$, and we shall write it $\dim
(\partial X)$.

\begin{corollary}
Let $Z$ be a nonempty $Q$-invariant analytic subset of $\tilde{X}^{\sigma }$
such that $\dim Z>\dim (\partial X)$; then $Z=\tilde{X}^{\sigma }$.
\end{corollary}

\begin{proof}
Let $d=\dim Z$, and let $Z^{(d)}$ be the set of $z\in Z$ such that $Z$ has a
component of dimension $d$ at $z$. Then $Z^{(d)}$ is an analytic subset of $Z
$ (Narasimhan 1966, p.~67) and is $Q$-invariant. Therefore, we can assume
that $Z=Z^{(d)}$. The image $Z^{\prime }$ of $Z$ in $X_{K}$ (any $K$) is
analytic and such that $Z^{\prime }=Z^{\prime (d)}$. As $\dim (X_{K}^{\ast
}\smallsetminus X_{K})<d$, the theorem of Remmert-Stein (loc. cit., p.~123)
shows that the closure $\bar{Z}^{\prime }$ of $Z$ in $X_{K}^{\ast }$ is
analytic, and Chow's theorem (loc. cit., p.~125) shows that $\bar{Z}^{\prime
}$ is algebraic. Therefore $\bar{Z}^{\prime }\cap X_{K}=Z^{\prime }$ is
algebraic, and the proposition applies.
\end{proof}

A point $\tilde{x}\in \tilde{X}$ is said to be \emph{special }if there
exists a torus $T\subset G$ such that $T_{\mathbb{C}}$ is maximal and $T(%
\mathbb{R})\tilde{x}=\tilde{x}$.

\begin{proposition}
Let $x\in \tilde{X}$ be special, and let $x^{\sigma }\in \tilde{X}^{\sigma }$
be any point such that $\hat{x}^{\sigma }=g(\sigma \hat{x})$ for some $g\in
G(\mathbb{A}_{f})$. Consider the $\sigma $-linear map of tangent spaces $%
\alpha \colon T_{x}(\tilde{X})\rightarrow T_{x^{\sigma }}(\tilde{X}^{\sigma
})$ defined by 
\begin{diagram}
x & \tilde{X} &  & \sigma \hat{X} & \rTo^{g} & \sigma \hat{X} & &
g^{-1}(x^{\sigma }) & \rMapsto & x^{\sigma } \\ 
\dMapsto & \dTo &  & \dTo &  &  & & \dMapsto &  &  \\ 
x_{K} & X_{K} & \rTo^{\sigma } & \sigma X_{K} &  &  & & \sigma x_{K} &  & 
\end{diagram}
Then there exists a homomorphism $j\colon T(\mathbb{Q})\rightarrow Q$ such
that $j(T(\mathbb{Q}{}))$ fixes $x$ and $\alpha $ commutes with the actions
of $T(\mathbb{Q}{})$ on the two tangent spaces. The closure of $j(T(\mathbb{Q%
}{}))$ in $\Aut(T_{x^{\sigma }}(\tilde{X}^{\sigma }))$ for the real topology
contains an element inducing multiplication by $\sqrt{-1}$.
\end{proposition}

\textsc{Proof.} Let $t\in T(\mathbb{Q}{})\subset G(\mathbb{A}_{f})$; then $t%
\hat{x}=\hat{x}$ and so $t\sigma \hat{x}=\sigma \hat{x}$. As $\sigma \hat{x}%
=g^{-1}\hat{x}^{\sigma }$, this means that $(gtg^{-1})\hat{x}^{\sigma }=\hat{%
x}^{\sigma }$. In particular, $gtg^{-1}$ maps one point of $\tilde{X}%
^{\sigma }$ into $\tilde{X}^{\sigma }$ and so stabilizes $\tilde{X}^{\sigma
} $. Thus $gtg^{-1}\in Q$, and we can define $j(t)=gtg^{-1}$. By
construction, $j(t)$ fixes $x^{\sigma }$, and it is routine to check that $%
\alpha $ carries the action of $t$ on $T_{x}(\tilde{X})$ into the action of $%
j(t)$ on $T_{x^{\sigma }}(\tilde{X}^{\sigma })$ ($\alpha =dg\circ
(dp_{K}^{\sigma })^{-1}\circ d\sigma \circ dp_{K}$; the actions of $t$ on $%
\tilde{X}$, $X_{K} $, and $\sigma X_{K}$ commute with $p_{K}$, $\sigma $,
and $p_{K}^{\sigma },$ and the action of $t$ on $\tilde{X}^{\sigma }$ is
transformed by $g\colon \tilde{X}^{\sigma }\rightarrow \tilde{X}^{\sigma }$
into the action of $gtg^{-1}$). To prove the last assertion of the
proposition, we need a lemma.

\begin{lemma}
Let $T$ be a torus over $\mathbb{Q}$. For each automorphism $\sigma $ of $%
\mathbb{C}$, there is a unique automorphism $t\mapsto t^{\sigma }$ of $T(%
\mathbb{C}{})$ such that $\chi (t^{\sigma })=\sigma \chi (t)$, all $t\in T(%
\mathbb{C})$, $\chi \in X^{\ast }(T)$. The map $t\mapsto t^{\sigma }$ is
continuous and, if $T(\mathbb{R}{})$ is compact, takes $T(\mathbb{R}{})$
into itself. For $t\in T(\mathbb{Q}{})$, $\sigma (\chi (t))=\chi (t^{\sigma
})$.
\end{lemma}

\begin{proof}
There is an isomorphism 
\begin{equation*}
t\mapsto (\chi \mapsto \chi (t))\colon T(\mathbb{C}{})\rightarrow \Hom%
(X^{\ast }(T),\mathbb{C}{}^{\times }),
\end{equation*}%
\noindent and we define $t\mapsto t^{\sigma }$ on $T(\mathbb{C}{})$ to
correspond to the map on $\Hom(X^{\ast }(T),\mathbb{C}{}^{\times })$ induced
by the $\mathbb{Z}{}$-linear map $\sigma \colon X^{\ast }(T)\rightarrow
X^{\ast }(T)$. Clearly, $\chi (t^{\sigma })=(\sigma \chi )(t)$, and $%
t\mapsto t^{\sigma }$ is the unique map with this property. The continuity
is obvious. If $T(\mathbb{R}{})$ is compact, it is the unique maximal
compact subgroup of $T(\mathbb{C}{})$, and so it is preserved by $t\mapsto
t^{\sigma }$. For $t\in T(\mathbb{Q}{})$, $\sigma \chi (t)=(\sigma \chi
)(t)=\chi (t^{\sigma })$.
\end{proof}

\textsc{Proof (of 1.10 continued): }The $\sigma $-linear map $\alpha \colon
T_{x}(\tilde{X})\rightarrow T_{x^{\sigma }}(\tilde{X}^{\sigma })$ induces a $%
\mathbb{C}{}$-linear isomorphism\footnote{%
For a complex vector space $V$, $\sigma V=V\otimes _{\mathbb{C}{},\sigma }%
\mathbb{C}{}$.} $\beta \colon \sigma T_{x}(\tilde{X})\rightarrow
T_{x^{\sigma }}(\tilde{X}^{\sigma })$. Moreover, 
\begin{equation*}
\beta \circ \sigma (\rho _{x}(t))=\rho _{x^{\sigma }}(j(t))\circ \beta
,\quad t\in T(\mathbb{Q}{}),
\end{equation*}
where $\rho _{x}$ and $\rho _{x^{\sigma }}$ denote the
representations of $G$ and $Q$ on the tangent spaces at $x$ and $x^{\sigma }$%
. Clearly, $(\sigma \rho _{x})(t)=\rho _{x}(t^{\sigma })$, and so 
\begin{equation*}
\beta \circ \rho _{x}(t^{\sigma })\circ \beta ^{-1}=\rho _{x^{\sigma
}}(j(t)),\quad t\in T(\mathbb{Q}{}).
\end{equation*}

\noindent As $T(\mathbb{Q}{})$ is dense in $T(\mathbb{R}{})$, it follows
that for any $\gamma \in T(\mathbb{R}{})$ there exists a $\gamma ^{\prime }$
in the closure of $\rho _{x^{\sigma }}(j(T(\mathbb{Q}{})))$ in $\Aut%
(T_{x^{\sigma }}(\tilde{X}^{\sigma }))$ such that $\gamma ^{\prime }$ acts
as $\beta \circ \rho _{x}({\gamma}^{\sigma })\circ \beta ^{-1}$ on $T_{x^{\sigma }}(%
\tilde{X}^{\sigma })$. It is known (see, for example, Helgason 1962, VIII
4.5) that there is a $\gamma $ in $T(\mathbb{R}{})$ acting as $\sqrt{-1}$ on 
$T_{x}(\tilde{X})$; therefore, 
\begin{eqnarray*}
\gamma ^{\prime } &=&\beta \circ \rho _{x}(\gamma^{\sigma} )\circ \beta ^{-1} \\
&=&\beta \circ (\text{multiplication by }\sqrt{-1})\circ \beta ^{-1} \\
&=&\text{multiplication by }\sqrt{-1}\text{.}\quad \quad \quad \square
\end{eqnarray*}

\section{A criterion to be an arithmetic variety.}

Let $(M,g)$ be an oriented Riemannian manifold. There is a unique volume
element $\mu $ on $M$ having value $1$ on any orthonormal frame. In local
coordinates 
\begin{eqnarray*}
g &=&\sum g_{ij}dx^{i}\otimes dx^{j} \\
\mu &=&\sqrt{\det (g_{ij})}dx^{1}\wedge dx^{2}\wedge \cdots \wedge dx^{m}.
\end{eqnarray*}

\noindent Let $\nabla $ be the connection defined by $g$. The \emph{%
curvature tensor }is defined by 
\begin{equation*}
R(X,Y)=[\nabla _{X},\nabla _{Y}]-\nabla _{\lbrack X,Y]}.
\end{equation*}

\noindent In terms of local coordinates, 
\begin{equation*}
R\left( \frac{\partial }{\partial x^{i}},\frac{\partial }{\partial x^{j}}%
\right) \frac{\partial }{\partial x^{\ell }}=\sum_{k}R_{\ell ij}^{k}\frac{%
\partial }{\partial x^{k}}.
\end{equation*}

\noindent The \emph{Ricci tensor }$r(X,Y)$ is determined in local
coordinates by 
\begin{equation*}
r(X,Y)=\sum R_{ij}dx^{i}\otimes dx^{j},\quad R_{ij}=\sum_{k}R_{ijk}^{k}.
\end{equation*}

Assume further that $M$ is a complex manifold, with multiplication by $\sqrt{%
-1}$ being described by the tensor $J$. Then $g$ is said to be \emph{%
Hermitian }if $g(JX,JY)=g(X,Y)$, all $X,Y$. In this case 
\begin{equation*}
h(X,Y)=_{\text{df}}g(X,Y)+ig(X,JY)
\end{equation*}

\noindent is a positive-definite Hermitian form. The form $\Phi
(X,Y)=_{\text{df}}g(X,JY) $ is skew-symmetric. It is called the fundamental $2$-form of
the Hermitian metric $h$, and if it is closed, $h$ is said to be \emph{K\"{a}%
hlerian}.

\begin{lemma}
Let $M$ be a connected complex manifold of dimension $n$ with K\"{a}hlerian
Riemannian structure $g$. If the associated volume element is described in
terms of local coordinates by 
\begin{equation*}
\mu =k\left( \frac{i}{2}\right) ^{n}dz_{1}\wedge d\bar{z}_{1}\wedge \cdots
\wedge d\bar{z}_{n}
\end{equation*}

\noindent then the Ricci tensor is described by 
\begin{equation}
r=\sum_{1\leq i,j\leq n}\frac{\partial ^{2}\log k}{\partial z_{i}\partial 
\bar{z}_{j}}dz_{i}\otimes d\bar{z}_{j}\quad \quad  \tag{2.1.1}
\end{equation}
\end{lemma}

\begin{proof}
Helgason 1962, VIII 2.5.
\end{proof}

In general, if $\mu =k(\frac{i}{2})^{n}dz_{1}\wedge d\bar{z}_{1}\wedge
\cdots \wedge d\bar{z}_{n}$ is a volume element on a complex manifold, then
the tensor defined by (2.1.1) will be called the \emph{Ricci tensor}, $%
Ric(\mu )$, of $\mu $.

\begin{lemma}
Let $M$ be a connected complex manifold on which a group $G$ acts
transitively, and let $\mu $ be a $G$-invariant volume element on $M$ such
that $g=Ric(\mu )$ is positive definite (and therefore is a Riemannian
metric). Then $g$ is equal to its Ricci tensor (and therefore is an Einstein
metric).
\end{lemma}

\begin{proof}
Let $\mu ^{\prime }$ be the volume form associated with $g$. Then $\mu
^{\prime }=f\mu $ for some positive function $f$ on $M$. As $\mu $ is $G$%
-invariant, so also are $g$, $\mu ^{\prime }$, and $f$. Therefore, $f$ is
constant, and so the Ricci tensor of $g$ is $Ric(f\mu )=Ric(\mu )=g$.
\end{proof}

\begin{theorem}
Let $M$ be a connected complex manifold on which a unimodular Lie group $G$
acts effectively and transitively. Assume there is a $G$-invariant volume
element $\mu $ on $M$ such that $Ric(\mu )$ is positive definite. Then $M$
is a Hermitian symmetric domain, and $G\supset \Aut(M)^{+}$.
\end{theorem}

\begin{proof}
This is proved in Koszul 1959, p.~61. Alternatively, one can combine the
following results:
\begin{quotation}
A complex homogeneous manifold with an invariant volume form whose Ricci
tensor is positive definite is isomorphic to a homogeneous bounded domain
(Piatetski-Shapiro 1969, p.~48).
\end{quotation}\begin{quotation}
A connected unimodular Lie group acting effectively and transitively on a
bounded domain is semisimple (Hano 1957).
\end{quotation}\begin{quotation}
A bounded domain admitting a transitive semisimple group of automorphisms is
symmetric (Borel 1954; Koszul 1955).
\end{quotation}
\end{proof}

\begin{theorem}
Let $M$ be a connected complex manifold, and endow $\Aut(M)$ with the
compact-open topology. Let $\Gamma \subset Q$ be subgroups of $\Aut(M)$ with 
$\Gamma $ discrete and torsion free, and assume that an orbit of $Q$ in $M$
is dense. Assume also that there is a $Q$-invariant volume form $\mu $ on $M$
such that

\begin{enumerate}
\item $Ric(\mu )$ is positive definite,

\item $\int_{\Gamma \backslash M}\mu <\infty $.
\end{enumerate}

\noindent Then $M$ is a Hermitian symmetric domain, and the closure $\bar{Q}$
of $Q$ in $\Aut(M)$ is a semisimple Lie group whose identity component is $%
\Aut(M)^{+}$. If moreover $Q$ is contained in the commensurability group of $%
\Gamma $ in $\bar{Q}$, then $\Gamma $ is arithmetic and so $\Gamma
\backslash M$ is an arithmetic variety.
\end{theorem}

\textsc{Proof. }Let $g=Ric(\mu )$. By assumption, $(M,g)$ is a Riemannian
manifold, and so the group $Is(M,g)$ of its isometries is a Lie group
(Kobayashi 1972, II Theorem 1.1). As $\mu $ is $Q$-invariant, $Q$ is
contained in $Is(M,g)$, and its closure $G$ is a Lie subgroup of $Is(M,g)$. Clearly, $G=\bar{Q}$.

\begin{lemma}
Let $(M,g)$ be a Riemannian manifold. Then, for all $m\in M$, the map $%
Is(M,g)\rightarrow M$, $\alpha \mapsto \alpha m$, is proper.
\end{lemma}

\begin{proof}
Let $O(M)$ be the bundle of orthonormal frames over $M$. Every automorphism $%
\alpha $ of $M$ defines a compatible automorphism $\bar{\alpha}$ of $O(M)$.
Let $u\in O(M)$ and let $m$ be its image in $M$. The mapping $\alpha \mapsto 
\bar{\alpha}(u)$ embeds $Is(M,g)$ as a closed submanifold of $O(M)$ (ibid.,
p.~41). The projection $O(M)\rightarrow M$ is clearly proper, and its
restriction to $Is(M,g)$ is the map $\alpha \mapsto \alpha m$.
\end{proof}

\textsc{Proof of 2.4 continued.} Let $m$ be a point of $M$ such that $Qm$ is
dense in $M$. Then the lemma shows that $Gm$ is closed and so equals $M$.
The orbit under $G^{+}$ of $m$ contains an open set in $M$, and so equals $M$
(Kobayashi and Nomizu 1963, Corollary 4.8, p.~178); thus $G^{+}$also acts
transitively on $M$. The isotropy group $K$ at the point $m$ is compact.
Thus the fact that $M=G^{+}/K$ carries a $G^{+}$-invariant measure such that 
$\Gamma \backslash M$ has finite volume implies that $G^{+}$ carries a left
invariant measure relative to which $\Gamma \cap G^{+}\backslash G^{+}$ has
finite volume. Hence $\Gamma \cap G^{+}$ is a lattice in $G^{+}$, and so $G$
is unimodular (the image of the modulus function $\Delta _{G}$ in $\mathbb{R}%
_{>0}$ is a subgroup with finite measure, and so is $\{1\}$ --- see
Raghunathan 1972, 1.9). The preceding theorem now shows that $M$ is a
Hermitian symmetric domain and $G^{+}=\Aut(M)^{+}$.

The final statement of the theorem is a consequence of the following theorem
of Margulis (Margulis 1977, Theorem 9):
\begin{quotation}
 Let $G$ be a semisimple connected real Lie group with no compact
factors, and let $\Gamma $ be a lattice in $G$; if the commensurability
group of $\Gamma $ in $G$ is dense in $G$, then $\Gamma $ is an arithmetic
subgroup of $G$ .\hspace*{\fill}$\square $
\end{quotation}
\begin{Remark}
In the statement of the theorem just applied, Margulis assumes that $\Gamma $
is irreducible, but this is unnecessary. If $\Gamma $ is reducible, then
there exist connected normal subgroups $G_{i}$ of $G$ such that $G=\prod
G_{i}$ (almost direct product), $\Gamma _{i}=_{\text{df}}\Gamma \cap G_{i}$
is an irreducible lattice in $G_{i}$, and $\prod \Gamma _{i}$ is a subgroup
of finite index in $\Gamma $. Clearly, the commensurability group of $\Gamma 
$ in $G$ is the product of the commensurability groups of the $\Gamma _{i}$: 
$\mathrm{Comm}(\Gamma )=\prod \mathrm{Comm}(\Gamma _{i})$. Thus, if $\mathrm{%
Comm}(\Gamma )$ is dense in $G$, each group $\mathrm{Comm}(\Gamma _{i})$ is
dense in $G_{i}$, and Margulis's statement can be applied to show that each $%
\Gamma _{i}$ is arithmetic.
\end{Remark}

\begin{corollary}
Let $X$ be an arithmetic variety as in \S 1, and suppose that there exists a 
$Q$-invariant volume element $\mu $ on $\tilde{X}^{\sigma }$ satisfying the
conditions (2.4a) and (2.4b). Assume that there exists a finite family of
arithmetic subvarieties $i_{\alpha }\colon X_{\alpha }\hookrightarrow X$
such that for some special point $x\in X$, $x_{\alpha }=_{\text{df}%
}i_{\alpha }^{-1}(x)$ is special in $X_{\alpha }$, all $\alpha $. If each $%
\sigma X_{\alpha }$ is arithmetic, and the subspace of $T_{x}(X)$ generated
by the $T_{x_{\alpha }}(X_{\alpha })$ has dimension $>\dim (\partial X)$,
then $\sigma X$ is arithmetic.
\end{corollary}

\begin{proof}
Let $g=Ric(\mu )$, and let $M$ be the closure of the orbit $Q\tilde{x}$ for
some $\tilde{x}\in \tilde{X}^{\sigma }$ lifting $x$. Then $(M,g)$ is a
Riemannian manifold, and the closure $\bar{Q}$ of $Q$ in $Is(M,g)$ acts
transitively on $M$ (by 2.5). Since $q\mapsto q\tilde{x}$, $\bar{Q}%
\rightarrow \tilde{X}^{\sigma }$, is a proper map, it is an embedding and $M$
is a regular closed submanifold of $\tilde{X}^{\sigma }$. By (1.10), $M$ has
a $\bar{Q}$-invariant complex structure --- it is therefore a complex
analytic subset of $\tilde{X}^{\sigma }$. Since it contains $\sigma \tilde{X}%
_{\alpha }$ for each $\alpha $, it has dimension $>\dim (\partial X)$, and
so (1.9) shows that $M=\tilde{X}^{\sigma }$.

Now (2.4) can be applied.
\end{proof}

\begin{remark}
It is now possible to complete the proof of Theorem 0.1 in the case that $X$
is compact. In this case, $\dim (\partial X)=-1$ and so the map $%
x\hookrightarrow X$ will serve for the family $(X_{\alpha }\hookrightarrow
X) $. Recall the following theorem of Yau (1978) (the Calabi conjecture):
\begin{quotation}
let $V$ be a smooth projective algebraic variety over $\mathbb{C}$ such that
$K_V$ is ample (equivalently, $c_1(V^{\text{an}})$ is negative); then there exists a unique K\"{a}hler metric on
$V^{\text{an}}$ such that $Ric(g)=g$.
\end{quotation}
Apply this to $\sigma X$ and let $\mu $ be the inverse image on $\tilde{X}^{\sigma }$ of the volume
element associated with $g$. The uniqueness of $g$ shows that $\mu $ is
invariant under all automorphisms of $\tilde{X}^{\sigma }$. Since $Ric(\mu
)=\tilde{g}$, it satisfies (a), and it satisfies (b) because $X$ is compact.
Thus (2.4) applies. \end{remark}

In the remainder of this section, we shall show that the condition (2.4b) is
automatically satisfied in the context of (2.6).

A complex manifold $M$ will be said to be \emph{compactifiable }if it can be
embedded as an open dense subset of a compact analytic space $\bar{M}$ in
such a way that $\bar{M}\smallsetminus M$ is an analytic subset. Hironaka's
theorem (Hironaka 1964) then shows that $\bar{M}$ can be chosen to be a manifold such
that $\bar{M}\smallsetminus M$ is a divisor with normal crossings. Clearly,
if $M=V^{\text{an}}$ with $V$ a quasi-projective algebraic variety, then $M$
is compactifiable (in fact, not even quasi-projectivity is necessary).

\begin{proposition}
Let $M$ be a compactifiable manifold of dimension $n$, and let $\mu $ be a
volume form on $M$ such that 
\begin{equation}
Ric(\mu )\text{ is positive-definite and }Ric(\mu )^{n}\geq \mu  \tag{2.8.1}
\end{equation}%
Then $\int_{M}\mu <\infty $.
\end{proposition}

Let 
\begin{eqnarray*}
D(a) &=&\{z\in \mathbb{C}\mid |z|<a\} \\
D^{\ast }(a) &=&\{z\in \mathbb{C}{}\mid 0<|z|<a\}.
\end{eqnarray*}

\begin{lemma}
For $\delta <1$,
\begin{eqnarray*}
\int_{D(\delta )}\frac{i}{\pi }\frac{dz\wedge d\bar{z}}{(1-|z|^{2})^{2}} &=&%
\frac{2\delta ^{2}}{1-\delta ^{2}} \\
\int_{D^{\ast }(\delta )}\frac{i}{\pi }\frac{dz\wedge d\bar{z}}{|z|^{2}(\log
|z|^{2})^{2}} &=&\frac{-1}{\log \delta }.
\end{eqnarray*}
\end{lemma}
\begin{proof}
Put $z=te^{i\theta }$, $\bar{z}=te^{-i\theta }$. Then $dz\wedge d\bar{z}%
=2itd\theta \wedge dt$. Thus, 
\begin{eqnarray*}
\int_{D(\delta )}\frac{i}{\pi }\frac{dz\wedge d\bar{z}}{(1-|z|^{2})^{2}}
&=&\int_{0}^{\delta }\frac{i}{\pi }2\pi \frac{2it}{(1-t^{2})^{2}}dt=\left. 
\frac{2}{1-t^{2}}\right| _{0}^{\delta }=\frac{2\delta ^{2}}{1-\delta ^{2}} \\
\int_{D^{\ast }(\delta )}\frac{i}{\pi }\frac{dz\wedge d\bar{z}}{|z|^{2}(\log
|z|^{2})^{2}} &=&\int_{0}^{\delta }\frac{i}{\pi }2\pi \frac{2it}{t^{2}4(\log
t)^{2}}dt=-\int \frac{dt}{t(\log t)^{2}}=\left. \frac{-1}{\log t}\right|
_{0}^{\delta }=\frac{-1}{\log \delta }.
\end{eqnarray*}
\end{proof}

\begin{remark}
Let 
\begin{equation*}
\mu _{D}=\frac{i}{\pi }\frac{dz\wedge d\bar{z}}{(1-|z|^{2})^{2}},\quad \mu
_{D^{\ast }}=\frac{i}{\pi }\frac{dz\wedge d\bar{z}}{|z|^{2}(\log |z|^{2})^{2}%
}.
\end{equation*}
Then $\mu _{D}$ is the Poincar\'{e} metric on the unit disc, and $Ric(\mu
_{D})=\mu _{D}$. The inverse image of $\mu _{D^{\ast }}$ relative to the
covering map $D\rightarrow D^{\ast }$ is $\mu _{D}$, and so $Ric(\mu
_{D^{\ast }})=\mu _{D^{\ast }}$ also. (Cf. Griffiths 1976, p.~47.)
\end{remark}

\begin{lemma}
Let $\mu $ be a volume element on $D(1)^{r}\times D^{\ast }(1)^{n-r}$
satisfying the estimates (2.8.1). Then 
\begin{equation*}
\mu \leq \mu _{D}^{r}\times \mu _{D^{\ast }}^{n-r}\text{.}
\end{equation*}
\end{lemma}

\begin{proof}
If $r=n$, this is precisely Ahlfors's lemma (Griffiths 1976, 2.21). Consider
the covering map 
\begin{equation*}
D(1)^{n}\overset{\varphi }{\rightarrow }D(1)^{r}\times D^{\ast }(1)^{n-r}%
\text{.}
\end{equation*}

\noindent Then $\varphi ^{\ast }(\mu )$ satisfies (2.8.1), and so $\varphi
^{\ast }(\mu )\leq \mu _{D}^{n}=\varphi ^{\ast }(\mu _{D}^{r}\times \mu
_{D^{\ast }}^{n-r})$. This implies that $\mu \leq \mu _{D}^{r}\times \mu
_{D^{\ast }}^{n-r}$.
\end{proof}

\begin{lemma}
Let $U=D(1)^{n}$ and let $U^{\prime }=D(1)^{r}\times D^{\ast }(1)^{n-r}$.
Then there exists an open neighbourhood $V$ of $U\smallsetminus U^{\prime }$
in $U$ such that $|\int_{V}\mu _{D}^{r}\times \mu _{D^{\ast }}^{n-r}|<\infty 
$.
\end{lemma}

\begin{proof}
A stronger statement is proved in (3.11) below.
\end{proof}

\begin{proof}[Proof of 2.8.]
Embed $M$ in a compact manifold $\bar{M}$ in such a way that $N=_{\text{df}}%
\bar{M}\smallsetminus M$ is a divisor with normal crossings. Then there is a
finite family of open subsets $U_{i}$ of $\bar{M}$ such that $N\subset \cup
U_{i}$ and, for each $i$, the pair $(U_{i},U_{i}\cap M)$ is isomorphic to $%
(D(1)^{n},D(1)^{r}\times D^{\ast }(1)^{n-r})$. For each $i$, choose $%
V_{i}\subset U_{i}$ to correspond to an open subset of $D(1)^{n}$ satisfying
the conditions of (2.12). Then the complement $C$ of $\cup V_{i}$ in $\bar{M}$ is compact, and
it is contained in $M$. Thus 
\begin{equation*}
\int_{M}\mu \leq \int_{C}\mu +\sum \int_{V_{i}\cap M}\mu \leq \int_{C}\mu
+\sum \int_{V_{i}\cap M}\mu _{D}^{r}\times \mu _{D^{\ast }}^{n-r}<\infty .
\end{equation*}
\end{proof}

\begin{corollary}
Let $X$ be an arithmetic variety (as in \S 1), and let $U$ be an open dense $%
Q$-invariant submanifold of $\tilde{X}^{\sigma }$ such that $\tilde{X}%
^{\sigma }\smallsetminus U$ is an analytic subset of $\tilde{X}^{\sigma }$.
Assume that there is a $Q$-invariant volume element $\mu$ on $U$ such that $%
Ric(\mu )>0$. Suppose further that there exists a finite family of
arithmetic subvarieties $i_{\alpha }\colon X_{\alpha }\hookrightarrow X$
such that, for some special point $x$ of $X$ such that $\sigma x$ is in the
image of $U\rightarrow \sigma X$, $i_{\alpha }^{-1}(x)$ is special in $%
X_{\alpha }$ for all $\alpha $. If each $\sigma X_{\alpha }$ is arithmetic,
and the subspaces $T_{x_{\alpha }}(X_{\alpha })$ of $T_{x}(X)$ generate it,
then $\sigma X$ is arithmetic.
\end{corollary}

\begin{proof}
Let $\tilde{x}\in U$ map to $\sigma x\in \sigma X$, and let $g=Ric(\mu )$
--- it is a Riemannian metric on $U$. Define $M$ to be the closure in $U$ of
the orbit $Q\tilde{x}$. As in the proof of (2.6), $M=\bar{Q}\tilde{x}$,
where $\bar{Q}$ is the closure of $Q$ in $Is(M,g)$, and it is a complex
analytic subset of $U$. As it contains the $X_{\alpha }$ for each $\alpha $,
it also open in $U$, and so $M=U$.

We know $\mu $ on $M$ satisfies condition (2.4a), and we now check that it
satisfies (2.4b). Suppose $Ric(\mu )^{n}=c\mu $ at some point of $\tilde{X}%
^{\sigma }$. Then $c>0$, and because $\bar{Q}$ acts transitively on $M$ and $%
\mu $ is $\bar{Q}$-invariant, we must have $Ric(\mu )^{n}=c\mu $ holding on
all of $M$. Since $Ric(c\mu )=Ric(\mu )$, this shows that $c\mu $ satisfies
the estimates (2.8.1). Therefore, so also does the volume element induced by 
$c\mu $ on $\Gamma ^{\sigma }\backslash M$, and so Proposition 2.8 implies
that $\int_{\Gamma ^{\sigma }\backslash M}c\mu <\infty $.

We can now apply (2.4) to $M$, $Q$, $\Gamma $, and we find that $\Gamma
\backslash M$ is an arithmetic variety. In particular, it is an open
algebraic subvariety of $\sigma X$, and so (1.8) implies that $M=\tilde{X}%
^{\sigma }$. We conclude that $\sigma X=\Gamma \backslash M$, and so it is
an arithmetic variety.
\end{proof}

\section{The Bergmann metric on $\tilde{X}^{\protect\sigma }$}

Let $M$ be a complex manifold of dimension $n$, and let
\begin{equation*}
\mathcal{H}(M)=\{w\in \Gamma (M,\Omega _{M}^{n\text{ hol}})\mid \left|
\tint_{M}\omega \wedge \bar{\omega}\right| <\infty \}
\end{equation*}

\noindent Then $\mathcal{H}(M)$ is a separable Hilbert space with inner
product 
\begin{equation*}
(\omega _{1}|\omega _{2})=i^{n^{2}}\int_{M}\omega _{1}\wedge \bar{\omega}_{2}%
\text{.}
\end{equation*}

\noindent Let $\omega _{0},\omega _{1},\ldots $ be an orthonormal basis for $%
\mathcal{H}(M)$, and let\footnote{Added 11.07.01:
For the convergence of $\sum \omega _{i}\wedge \bar{\omega}_{i}$, see, for example,
Weil, A., Introduction \`a l'\'etude des vari\'et\'es k\"ahl\'eriennes,
 Hermann, Paris, 1958, Chapter III, Thm 1, p.~60.}
\begin{equation*}
\mu _{M}=\sum \omega _{i}\wedge \bar{\omega}_{i}.
\end{equation*}

\noindent This is a nonnegative $C^{\infty }$ $2n$-form on $M$, called the 
\emph{Bergmann volume form}. When it has no zeros it is a volume element.
Note that 
\begin{equation*}
\int_{M}\mu _{M}=\dim \mathcal{H}{}(M).
\end{equation*}

\begin{proposition}
(a). Let $m\in M$; then 
\begin{equation*}
\mu _{M}(m)=\sup_{\substack{ \omega \in \mathcal{H}{}(M)  \\ (\omega |\omega
)=1}}(\omega \wedge \bar{\omega})(m).
\end{equation*}

(b) All automorphisms of $M$ leave $\mu _{M}$ invariant.

(c) If $M^{\prime }$ is a connected open submanifold of $M$, then $\mu
_{M}|M^{\prime }=c\mu _{M^{\prime }}$ where $c$ is a function on $M^{\prime
} $, $0\leq c\leq 1$; moreover, if $M\smallsetminus M^{\prime }$ is a
complex analytic subvariety of dimension $\leq n-1$, then $\mathcal{H}%
{}(M)\rightarrow \mathcal{H}{}(M^{\prime })$ is bijective and so $\mu
_{M}|M^{\prime }=\mu _{M^{\prime }}$.

(d) If $M_{1}$ and $M_{2}$ are complex manifolds of dimensions $n_{1}$ and $%
n_{2}$, then 
\begin{equation*}
\mu _{M_{1}\times M_{2}}=(-1)^{n_{1}n_{2}}\mu _{M_{1}}\wedge \mu _{M_{2}}.
\end{equation*}
\end{proposition}

\begin{proof}
Kobayashi 1959, 2.1, 2.2, 2.3, 2.4, 2.5.
\end{proof}

Consider the condition:

\begin{quotation}
(3.2) for every $m\in M$, there exists an $\omega \in \mathcal{H}{}(M)$ such
that $\omega (m)\neq 0$.
\end{quotation}

\noindent When this condition holds, $\mu _{M}$ is a volume element, and we
let $h_{M}=Ric(\mu _{M})$ be the associated Hermitian tensor.

Consider the condition:

\begin{quotation}
(3.3) for all $m\in M$ and $Z\in T_{m}(M)$, there exists an $\omega
=fdz_{1}\wedge \cdots \wedge dz_{m}$ in $\mathcal{H}{}(M)$ such that $f(m)=0$%
, $Z(f)\neq 0$.
\end{quotation}
\setcounter{X}{3}

\begin{proposition}
The form $h_{M}$ is invariant under all automorphisms of $M$ and is positive
semi-definite. It is positive definite if and only if (3.3) holds.
\end{proposition}
\begin{proof}
Ibid.~3.1.
\end{proof}

\begin{remark}
Let $\mathbb{P}{}(\mathcal{H}{}(M)^{\vee })$ be the (possibly infinite
dimensional) projective space of lines in the dual Hilbert space to $%
\mathcal{H}{}(M)$, and assume that $M$ satisfies (3.2). Then there is a
canonical map $j\colon M\rightarrow \mathbb{P}{}(\mathcal{H}{}(M)^{\vee })$
such that, if $\omega _{0}$ is nowhere zero on $U\subset M$, then $j(m)$, $%
m\in U$, is the class of the map $\omega \mapsto (\omega /\omega _{0})(m)$.
It is possible to regard $\mathbb{P}{}(\mathcal{H}{}(M)^{\vee })$ as an
infinite-dimensional complex manifold, and the usual construction in the
finite-dimensional case generalizes to give a complete K\"{a}hler metric on $%
\mathbb{P}{}(\mathcal{H}{}(M)^{\vee })$. The Bergmann Hermitian form $h_{M}$
on $M$ is the inverse image under $j$ of the canonical metric on $\mathbb{P}%
{}(\mathcal{H}{}(M)^{\vee })$. The map $j$ is an immersion, and $h_{M}$ is a
metric, if and only if (3.3) holds. (See Kobayashi 1959, \S 7,8 for this.)
\end{remark}

\begin{remark}
Let $X$ be an arithmetic variety, and assume that there exists a family $%
(X_{\alpha }\hookrightarrow X)$ as in (2.6). Then it follows from (2.6) and
(2.12) that $\sigma X$ is an arithmetic variety if the Bergmann Hermitian
form on $\tilde{X}^{\sigma }$ is a metric; conversely, if $\sigma X$ is an
arithmetic variety, then $\tilde{X}^{\sigma }$ is a bounded domain and so
Bergmann's original theorem says that it has a nondegenerate Bergmann metric.
\end{remark}

The main result of this section is a first step toward proving that $\tilde{X%
}^{\sigma }$ has a nondegenerate Bergmann metric.

\begin{theorem}
With the notations of \S 1, $\mathcal{H}{}(\tilde{X}^{\sigma })\neq 0$.
\end{theorem}

\begin{proposition}
Let $M$ be a compactifiable complex manifold. For all $\varepsilon >0$,
there exists an open subset $U_{\varepsilon }$ in $M$ with compact
complement such that, for any \'{e}tale covering $\varphi \colon
N\rightarrow M$, 
\begin{equation*}
\left| \int_{\varphi ^{-1}(U_{\varepsilon })}\mu _{N}\right| <\varepsilon
\deg \varphi
\end{equation*}

\noindent where $\mu _{N}$ is the Bergmann volume form on $N$.
\end{proposition}

The proof is based on the following elementary result.

\begin{lemma}
Let $\varphi _{m}\colon D^{\ast }(1)\rightarrow D^{\ast }(1)$ be the map $%
z\mapsto z^{m}$. For all $\varepsilon >0$, there exists a $\delta $ such
that 
\begin{equation*}
\left| \int_{\varphi _{m}^{-1}(D^{\ast }(\delta ))}\mu _{D^{\ast
}(1)}\right| <\varepsilon m,\quad \text{all }m.
\end{equation*}
\end{lemma}

\begin{proof}
The Bergmann metric on $D(1)$ is the Poincar\'{e} metric, and so from (3.1c)
and (2.9) we find that 
\begin{equation*}
\int_{\varphi _{m}^{-1}(D^{\ast }(\delta ))}\mu _{D^{\ast
}(1)}=\int_{D^{\ast }(\delta ^{\frac{1}{m}})}\frac{i}{\pi }\frac{dz\wedge d%
\bar{z}}{(1-|z|^{2})^{2}}=\frac{2}{(1/\delta )^{2/m}-1}.
\end{equation*}

As $x-1>\log x$ for all $x>1$, we see that 
\begin{equation*}
\frac{2}{(1/\delta )^{2/m}-1}<\frac{2}{\log (1/\delta )^{2/m}}=\frac{m}{\log
(1/\delta )}
\end{equation*}

from which the lemma is obvious.
\end{proof}

\begin{remark}
Note that the notation in this section conflicts with that in (2.10) --- $%
\frac{i}{\pi }\frac{dz\wedge d\bar{z}}{|z|^{2}(\log |z|^{2})^{2}}$ is not
the Bergmann volume element on $D^{\ast }(1)$. Denote $\frac{i}{\pi }\frac{%
dz\wedge d\bar{z}}{|z|^{2}(\log |z|^{2})^{2}}$ by $\nu _{D^{\ast }(1)}$.
Then Ahlfors's lemma applied on the covering space $D(1)$ of $D^{\ast }(1)$
shows that $\nu _{D^{\ast }(1)}\geq \mu _{D^{\ast }(1)}$ (cf. 2.11). Thus,
we find again that 
\begin{equation*}
\int_{D^{\ast }(\delta ^{1/m})}\mu _{D^{\ast }(1)}\leq \int_{D^{\ast
}(\delta ^{1/m})}\nu _{D^{\ast }(1)}\overset{(2.9)}{=}\frac{1}{\log
(1/\delta )^{m}}=\frac{m}{\log (1/\delta )}\text{.}
\end{equation*}
\end{remark}

\begin{lemma}
Let $U=D(1)^{m}$ and $U^{\prime }=D(1)^{r}\times D^{\ast }(1)^{n-r}$. For
all $\varepsilon >0$, there exists an open neighbourhood $V$ of $%
U\smallsetminus U^{\prime }$ on $U$ such that for any \'{e}tale
covering\footnote{That is, a finite \'etale covering.} $\varphi \colon N\rightarrow U^{\prime }$, $|\int_{\varphi
^{-1}(V)}\mu _{N}|<\varepsilon \deg \varphi$.
\end{lemma}

\begin{proof}
The fundamental group of $U^{\prime }$ is $\mathbb{Z}{}^{n-r}$, and so every 
\'{e}tale covering of it is of the form%
\begin{equation*}
N=D(1)^{r}\times D^{\ast }(1)^{n-r}\rightarrow D(1)^{r}\times D^{\ast
}(1)^{n-r},\quad (z_{1},\ldots )\mapsto (z_{1},\ldots
,z_{r},z_{r+1}^{m_{1}},\ldots ).
\end{equation*}

\noindent The proof is too messy\footnote{Added 11.07.01: Chai suggested (in 1988) that
``a short proof of the lemma seems possible''.} to write out in detail. Consider for
example the case $n=2$, $U^{\prime }=D^{\ast }(1)^{2}$. Then 
\begin{equation*}
\text{\lbrack Picture omitted]}
\end{equation*}

\noindent (a) Look at $V\cap D(\delta )\times D(\delta )$ --- apply (3.9).

\noindent (b) Each other connected component $V_{0}$ of $V$ is simply connected. Thus 
\begin{equation*}
\varphi ^{-1}(V_{0})=\text{disjoint union of }\deg \varphi \text{ copies of }%
V_{0}.
\end{equation*}

Hence 
\begin{equation*}
\int_{\varphi ^{-1}(V_{0})}\mu _{N}\leq \int_{\varphi ^{-1}(V_{0})}\mu
_{\varphi ^{-1}(V_{0})}=\deg (\varphi )\int_{V_{0}}\mu _{V_{0}}\text{,}
\end{equation*}

\noindent and so one only has to arrange things so that 
\begin{equation*}
\int_{V_{0}}\mu _{V_{0}}<\frac{\varepsilon }{100}
\end{equation*}
\end{proof}

\begin{proof}[Proof of 3.8] Embed $M$ in a compact manifold $\bar{M}$ in such a way that $\bar{%
M}\smallsetminus M$ is a divisor with normal crossings. Then there is a
finite family of open subsets $(U_{i})_{1\leq i\leq s}$ of $\bar{M}$ such
that $\bar{M}\smallsetminus M\subset \cup U_{i}$ and, for each $i$, the pair 
$(U_{i},U_{i}\cap M)$ is isomorphic to $(D(1)^{n},D(1)^{r}\times D^{\ast
}(1)^{n-r})$. For each $i$, choose a neighbourhood $V_{i}$ of $%
U_{i}\smallsetminus U_{i}\cap M$, as in the sublemma, for $\varepsilon /s$.
Then the complement of $U=_{\text{df}}\cup V_{i}$ on $\bar{M}$ is compact,
and it is contained in $M$. Moreover, for any \'{e}tale covering $\varphi
\colon N\rightarrow M$, 
\begin{equation*}
\left| \int_{\varphi ^{-1}(U)}\mu _{N}\right| \leq \sum_{i=1}^{s}\left|
\int_{\varphi ^{-1}(V_{i})}\mu _{N}\right| <s(\frac{\varepsilon }{s}\deg
\varphi )=\varepsilon \deg \varphi .
\end{equation*}
\end{proof}

The next result gives a criterion for showing $\mathcal{H}{}(\tilde{M})$ is
nonzero.

\begin{proposition}
Let $M$ be a compactifiable manifold of dimension $n$, and let $p\colon 
\tilde{M}\rightarrow M$ be an infinite Galois covering with Galois group $%
\Gamma $. Assume there is a sequence of normal subgroups of finite index in $%
\Gamma $, 
\begin{equation*}
\Gamma =\Gamma _{0}\supset \Gamma _{1}\supset \cdots \supset \Gamma
_{i}\supset \cdots \supset \{1\}
\end{equation*}

\noindent such that $\cap \Gamma _{i}=\{1\}$.

Let $M_{i}=\Gamma _{i}\backslash \tilde{M}$ and let $h_{i}=\dim \mathcal{H}%
{}(M_{i})$. Then

\begin{enumerate}
\item $h_{i}<\infty $,

\item $\{h_{i}/(\Gamma \colon \Gamma _{i})\}$ is bounded,

\item if the sequence $h_{i}/(\Gamma \colon \Gamma _{i})$ does not tend to
zero, then $\mathcal{H}{}(\tilde{M})\neq 0$.
\end{enumerate}
\end{proposition}

\begin{remark}
In general, if $\varphi \colon Y\rightarrow X$ is an \'{e}tale covering, one
can show that 
\begin{equation*}
(\deg \varphi )\varphi ^{\ast }(\mu _{X})\geq \mu _{Y}\geq \left( \frac{1}{%
\deg \varphi }\right) \varphi ^{\ast }(\mu _{X}).
\end{equation*}%
\noindent \noindent Hence 
\begin{equation*}
(\deg \varphi )^{2}\int_{X}\mu _{X}\geq \int_{Y}\mu _{Y}\geq \int_{X}\mu
_{X}.
\end{equation*}%
\noindent \noindent For example, if $Y$ is a disjoint union of copies of $X$%
, then $\mu _{Y}=\varphi ^{\ast }(\mu _{X})$ and $\int_{Y}\mu _{Y}=\deg
\varphi \cdot \int_{X}\mu _{X}$.

\noindent In the situation of the proposition, one has trivial estimates 
\begin{equation*}
(\Gamma \colon \Gamma _{i})h_{0}\geq \frac{h_{i}}{(\Gamma \colon \Gamma _{i})%
}\geq \frac{h_{0}}{(\Gamma \colon \Gamma _{i})}
\end{equation*}%
\noindent \noindent with $h_{i}=(\Gamma \colon \Gamma _{i})$ in the case
that $\tilde{M}$ is a trivial covering of $M$.
\end{remark}

\begin{proof}[Proof of 3.12] As the Bergmann volume form $\mu _{M_{i}}$ on $M_{i}$ is $\Gamma $%
-invariant, it induces a form $\mu _{i}$ on $M$ such that $q_{i}^{\ast }\mu
_{i}=\mu _{M_{i}}$ where $q_{i}$ is the covering map $M_{i}\rightarrow M$.
Clearly, 
\begin{equation}
\int_{M}\mu _{i}=\frac{1}{(\Gamma \colon \Gamma _{i})}\int_{M_{i}}\mu
_{M_{i}}=\frac{h_{i}}{(\Gamma \colon \Gamma _{i})}  \tag{3.12.1}
\end{equation}

\noindent According to (3.8), there exists an open subset $U=U_{1}$ with
compact complement, $C=M\smallsetminus U$, such that 
\begin{equation*}
\left| \int_{\varphi ^{-1}(U)}\mu _{N}\right| <\deg \varphi
\end{equation*}

\noindent for any \'{e}tale covering $\varphi \colon N\rightarrow M$.
Because $C$ is compact, there exists a finite set $(U_{r})_{1\leq r\leq s}$
of open subsets of $M$ such that

\begin{enumerate}
\item there exists isomorphisms $\varphi _{r}\colon U_{r}\overset{\approx }{%
\rightarrow }D(1)^{n}$ (of complex manifolds);

\item $\bigcup_{r=1}^{s}\varphi _{r}^{-1}(D(\frac{1}{2})^{n})\supset C$.
\end{enumerate}

Because $q_{i}^{-1}(U_{r})$ is a disjoint union of copies of $U_{r}$, we
have 
\begin{equation*}
q_{i}^{\ast }(\mu _{U_{r}})=\mu _{q_{i}^{-1}(U_{r})}\geq \mu _{M_{i}}
\end{equation*}

\noindent and so $\mu _{U_{r}}\geq \mu _{i}|U_{r}$. Note that $\mu
_{U_{r}}=\varphi _{r}^{\ast }(\mu _{D(1)}^{n})$ (see 3.1d). We conclude that 
\begin{equation*}
\int_{C}\mu _{i}\leq s\int_{D(\frac{1}{2})^{n}}\mu _{D(1)^{n}}=s\left(
\int_{D(\frac{1}{2})}\mu _{D(1)}\right) ^{n}=B<\infty ,
\end{equation*}

\noindent and that 
\begin{equation*}
\frac{h_{i}}{(\Gamma \colon \Gamma _{i})}\overset{3.12.1}{=}\int_{M}\mu
_{i}=\int_{C}\mu _{i}+\int_{U}\mu _{i}\leq B+1\text{.}
\end{equation*}

\noindent This proves both (a) and (b).

Now assume that the sequence $h_{i}/(\Gamma \colon \Gamma _{i})$ does not
tend to zero. By passing to a subsequence, we can in fact assume that for
some $a>0$, $h_{i}/(\Gamma \colon \Gamma _{i})\geq a$ all $i$. Choose $%
U\subset M$ as in (3.8) with $\varepsilon =a/2$, and let $C=M\smallsetminus
U $; then 
\begin{equation*}
\int_{C}\mu _{i}=\int_{M}\mu _{i}-\int_{U}\mu _{i}=\frac{h_{i}}{(\Gamma
\colon \Gamma _{i})}-\int_{\varphi ^{-1}(U)}\mu _{M_{i}}\geq a-\frac{a}{2}=%
\frac{a}{2}\text{.}
\end{equation*}

\noindent Let $\nu =\sum_{r=1}^{s}\varphi _{r}^{\ast }(\mu _{D(1)^{n}})$. We
showed above that, on $U_{r}$, 
\begin{equation*}
\mu _{i}\leq \mu _{U_{r}}=\varphi _{r}^{\ast }(\mu _{D(1)^{n}}),
\end{equation*}

\noindent and so, on $C$, $\mu _{i}\leq \nu $.
\end{proof}

\begin{lemma}
There exists an $x_{0}\in C$ such that a subsequence of $(\mu _{i}/\nu
)(x_{0})$ converges to some $b>0$.
\end{lemma}

\begin{proof}
Suppose $\sup_{x\in C}((\mu _{i}/\nu )(x))\rightarrow 0$; then $\int_{C}\mu
_{i}\rightarrow 0$, which contradicts the assertion above that $\int_{C}\mu
_{i}=a/2$. Hence the sequence $\sup ((\mu _{i}/\nu )(x))$ does not tend to
zero, and so we can choose a subsequence of $i$'s for which $\sup_{x}((\mu
_{i}/\nu )(x))\rightarrow b>0$. Choose $x_{i}$ such that $|(\mu _{i}/\nu
)(x_{i})-\sup ((\mu _{i}/\nu )(x))|<1/i$. Then $(\mu _{i}/\nu
)(x_{i})\rightarrow b$. Now take a subsequence of the $x_{i}$ converging to
a limit $x_{0}$. Then clearly, $(\mu _{i}/\nu )(x_{0})\rightarrow b$.
\end{proof}

\begin{lemma}
Let $\tilde{x}_{0}\in \tilde{M}$. Then there exists a sequence of open
neighbourhoods 
\begin{equation*}
\cdots \subset N_{i}\subset N_{i+1}\subset \cdots
\end{equation*}

\noindent of $\tilde{x}_{0}$ such that $\tilde{M}=\cup N_{i}$ and the
restriction of $p_{i}\colon \tilde{M}\rightarrow M_{i}$ to $N_{i}$ is an
isomorphism of $N_{i}$ with a dense open subset of $M_{i}$.
\end{lemma}

\begin{proof}
Choose a Riemannian metric $\rho $ on $M$, and let $\tilde{\rho}$ be the
metric it induces on $\tilde{X}^{\sigma }$. Define 
\begin{equation*}
N_{i}=\{\tilde{x}\in \tilde{M}\mid \tilde{\rho}(\tilde{x}_{0},\tilde{x})<%
\tilde{\rho}(\tilde{x}_{0},\tilde{x})\text{, all }\gamma \in \Gamma _{i}%
\text{, }\gamma \neq 1\}.
\end{equation*}

\noindent Check that these sets have the right property (cf. Kazhdan 1975, p.~167).
\end{proof}

We now prove (c) of the proposition. Let $x_{0}$ be as in (3.14), and let $%
\tilde{x}_{0}\in \tilde{M}$ map to it. Then, from (3.1a) we know that there
exists an $\omega _{i}\in \mathcal{H}{}(M_{i})$ such that $(\omega
_{i}|\omega _{i})=1$ and 
\begin{equation}
(\omega _{i}\wedge \bar{\omega}_{i})(p_{i}(\tilde{x}_{0}))\geq \frac{1}{2}%
\mu _{M_{i}}(p_{i}(\tilde{x}_{0}))\text{.}  \tag{3.15.1}
\end{equation}

\begin{lemma}
There exists a compact neighbourhood $\tilde{U}$ of $\tilde{x}_{0}$ such
that for some $c>0$, $\int_{p_{i}(\tilde{U})}\omega _{i}\wedge \bar{\omega}%
_{i}>c$, all $i$.
\end{lemma}

\begin{proof}
Obvious from (3.14) and (3.15.1). (Cf. Kazhdan 1983, p.~153.).
\end{proof}

Now let $\hat{\mathcal{H}}$ be the Hilbert space of measurable
square-integrable sections $\eta $ of $\Omega _{\tilde{M}}^{n}$. For each $i$%
, define $\eta _{i}\in \hat{\mathcal{H}}$ by 
\begin{equation*}
\left\{ 
\begin{array}{lll}
\eta _{i}|N_{i} & = & p_{i}^{\ast }(\omega _{i}) \\ 
\eta _{i}|(\tilde{M}\smallsetminus N_{i}) & = & 0.%
\end{array}%
\right.
\end{equation*}

\noindent As $(\eta _{i}|\eta _{i})=1$, there exists a weakly convergent
subsequence of the $\eta _{i}$ tending to $\eta \in \hat{\mathcal{H}}$. Now $%
\eta $ is holomorphic on $\cup N_{i}=\tilde{M}$, and $\int_{\tilde{U}}\eta
\wedge \bar{\eta}>c/2$, and so $\eta \neq 0$. Thus $\mathcal{H}{}(\tilde{M}%
)\neq 0$.

\begin{conjecture}
Let $\bar{G}$ be a semisimple real Lie group, $\Gamma $ an arithmetic
subgroup of $\bar{G}$, and 
\begin{equation*}
\Gamma =\Gamma _{0}\supset \Gamma _{1}\supset \cdots \supset \Gamma
_{i}\supset \cdots \supset \{1\}
\end{equation*}

\noindent a sequence of normal subgroups of $\Gamma $ of finite index such
that $\cap \Gamma _{i}=\{1\}$. Let $W$ be an irreducible cuspidal
representation of $\bar{G}$ and define $h_{i}(W)=\dim (\Hom_{\bar{G}%
}(W,L^{2}(\Gamma _{i}\backslash \bar{G}))$. Then $h_{i}(W)/(\Gamma \colon
\Gamma _{i})$ does not tend to zero as $i\rightarrow \infty $.
\end{conjecture}

In Kazhdan 1983, p.~156, Kazhdan calls this ``Theorem A, which we will prove
in another paper''. At present the statement still seems to be unproven, but
Clozel has recently shown the following weaker result which is sufficient
for our purposes:

\begin{quotation}
Let $G$ be a simply connected semisimple algebraic group over $\mathbb{Q}{}$%
, and let $p_{0}$ be a prime such that $G(\mathbb{Q}{}_{p_{0}})$ has a
supercuspidal representation $\pi _{p_{0}}$ (e.g., take $p_{0}$ to be any
prime such that $G$ is split over $\mathbb{Q}{}_{p_{0}}$); let $K_{p_{0}}$be
a compact open subgroup of $G(\mathbb{Q}{}_{p_{0}})$ such that $\mathcal{H}%
{}(\pi _{p_{0}})^{K_{p_{0}}}\neq 0$ (here $\mathcal{H}{}$ is the Hecke
algebra); let $K_{0}\supset K_{1}\supset \cdots $ be a sequence of compact
open subgroups of $G(\mathbb{A}_{f})$ such that

(a) $(K_{i})_{p_{0}}=K_{p_{0}}$;

(b) there exists a finite set $S$ of primes such that $(K_{i})_{p}$ is
maximal for $p\notin S$;

(b) $\cap K_{i}=\{1\}$.

\noindent Then, for any irreducible cuspidal representation $W$ of $G(%
\mathbb{R}{})$, there exists a constant $a>0$ such that 
\begin{equation*}
\frac{h_{i}(W)}{(\Gamma _{0}\colon \Gamma _{i})}\geq a,\quad \Gamma _{i}=G(%
\mathbb{Q}{})\cap K_{i}.
\end{equation*}
\end{quotation}

Now return to the situation in the statement of (3.7). Then $W=_{\text{df}}%
\mathcal{H}{}(\tilde{X})$ is in a natural way a cuspidal representation of $%
G(\mathbb{R}{})$.

\begin{lemma}
There is an isomorphism 
\begin{equation*}
\Hom_{G(\mathbb{R}{})}(W,L^{2}(\Gamma \backslash \tilde{X}))\overset{\approx 
}{\rightarrow }\mathcal{H}{}(X).
\end{equation*}
\end{lemma}

\begin{proof}
Well known --- see Kazhdan 1983, p.~156--157.
\end{proof}

\begin{lemma}
The dimensions of $\mathcal{H}{}(X)$ and $\mathcal{H}{}(\sigma X)$ are equal.
\end{lemma}

\begin{proof}
Let $\bar{X}$ be a smooth variety containing $X$ as a dense open subvariety.
Then (3.1c) shows that $\mathcal{H}{}(X)\overset{\approx }{\rightarrow }%
\mathcal{H}{}(\bar{X})$, and $\mathcal{H}{}(\bar{X})=\Gamma (\bar{X},\Omega
_{X,\text{alg}}^{n})$. The lemma is now obvious.
\end{proof}

We are now ready to prove (3.7). Choose a family $K_{i}$ satisfying the
conditions of Clozel's theorem. Then the theorem and (3.19) show that for
some $a>0$, $\dim \mathcal{H}{}(X_{i})/(\Gamma _{0}\colon \Gamma _{i})\geq a$
all $i$, where $X_{i}=\Gamma (K_{i})\backslash \tilde{X}$. Now (3.20) shows
that $\dim \mathcal{H}{}(\sigma X_{i})/(\Gamma _{0}^{\sigma }\colon \Gamma
_{i}^{\sigma })\geq a$, where $\Gamma _{i}^{\sigma}=\Gamma ^{\sigma }(K_{i})$, and
(3.12) implies that $\mathcal{H}{}(\tilde{X}^{\sigma })\neq 0$.

\section{Subbundles of $T(\tilde{X})$}

Let $G$ be a simply connected semisimple algebraic group over $\mathbb{Q}{}$%
, and let $\tilde{X}$ be a symmetric Hermitian domain on which $G(\mathbb{R}%
{})$ acts in such a way that 
\begin{equation*}
\tilde{X}\approx G(\mathbb{R}{})/\{\text{maximal compact subgroup}\}.
\end{equation*}

\noindent Let 
\begin{equation*}
\tilde{X}=X_{1}\times \cdots \times X_{s},\quad X_{i}\text{ irreducible
symmetric Hermitian domain.}
\end{equation*}

\noindent Choose a point $x=(x_{1},\ldots ,x_{s})\in \tilde{X}$. Then 
\begin{eqnarray*}
G_{\mathbb{R}{}} &=&G_{1}\times \cdots \times G_{s}\times K_{0},\quad G_{i}%
\text{ noncompact, }K_{0}\text{ compact,} \\
G_{i}(\mathbb{R}{})x &=&G_{i}(\mathbb{R})x_{i}=X_{i}\text{.}
\end{eqnarray*}

\noindent Let 
\begin{equation*}
X_{i}^{\prime }=\{x_{1}\}\times \cdots \times X_{i}\times \cdots \times
\{x_{s}\}\subset \tilde{X}.
\end{equation*}

\noindent Then $T(X_{i}^{\prime })\subset T(\tilde{X})$, and $T(\tilde{X}%
)=\oplus _{i=1}^{s}T(X_{i}^{\prime })$. For $I\subset \{1,\ldots ,s\}$,
define 
\begin{equation*}
T^{I}(\tilde{X})=\oplus _{i\in I}T(X_{i}^{\prime })=\{(x,v)\in T(\tilde{X}%
)\mid v\in \oplus T_{x_{i}}(X_{i})\}.
\end{equation*}

\noindent Then $T^{I}(\tilde{X})$ is a subbundle of $T(\tilde{X})$, stable
under $G(\mathbb{R}{})$.

\begin{lemma}
Every $G(\mathbb{Q}{})$-stable complex subbundle of $T(\tilde{X})$ is of the
form $T^{I}(\tilde{X})$ for some $I\subset \{1,\ldots ,s\}$.
\end{lemma}

\begin{proof}
Let $W$ be such a subbundle. As $G(\mathbb{Q}{})$ is dense in $G(\mathbb{R}%
{})$ (real approximation theorem) and $G(\mathbb{R}{})$ acts transitively on 
$\tilde{X}$, it will suffice to show that $W_{x}=T^{I}(\tilde{X})_{x}$, some 
$I$, for a fixed $x\in \tilde{X}$. Let $K_{x}$ be the isotropy group at $x$.
Then $W_{x}$ is a $K_{x}$-stable subspace of $T(\tilde{X})_{x}=\oplus
T(X_{i})_{x_{i}}$. With the usual notations, let 
\begin{eqnarray*}
\Lie(G) &=&\mathfrak{g}=\mathfrak{k}{}+\mathfrak{p}{},\quad \mathfrak{k}{}=%
\Lie(K_{x}),\quad \mathfrak{p}{}=T_{x}(\tilde{X}), \\
\Lie(G_{i}) &=&\mathfrak{g}_{i}=\mathfrak{k}{}{}_{i}+\mathfrak{p}%
{}_{i},\quad \mathfrak{p}_{i}=T_{x_{i}}(\tilde{X}_{i}).
\end{eqnarray*}
Then 
\begin{equation*}
\mathfrak{k}{}=\Lie(K_{0})\oplus \mathfrak{k}{}_{1}\oplus \cdots \oplus 
\mathfrak{k}{}_{s}
\end{equation*}
and 
\begin{equation*}
\mathfrak{p}{}=\oplus \mathfrak{p}{}_{i}.
\end{equation*}%
Note that $W_{x}$ is a $\mathfrak{k}{}$-stable subspace of $\mathfrak{p}{}$.
Almost by definition of what it means for $X_{i}$ to be irreducible
(Helgason 1962, VIII.5, p.~377), the action of $\mathfrak{k}{}_{i}$ on $\mathfrak{p%
}{}_{i}$ is irreducible. Thus $W_{x}=\oplus _{i\in I}\mathfrak{p}_{i}$, some $I$.
\end{proof}

\begin{lemma}
Let $G$ be a simply connected, almost simple, algebraic group over a number
field $F$. Let $S_{\infty }$ be the set of infinite primes of $F$, and let $%
v\in S_{\infty }$ be such that $G(F_{v})$ is not compact. Then for any
congruence group $\Gamma \subset G(F)$, $G(F_{v})\Gamma $ is dense in $%
\prod_{v\in S_{\infty }}G(F_{v})$.
\end{lemma}

\begin{proof}
For some compact open subgroup $K$ of $G(\mathbb{A}_{f,F})$, $\Gamma \supset
K\cap G(F)$. The strong approximation theorem shows that $G(F_{v})G(F)$ is
dense in 
\begin{equation*}
G(\mathbb{A}{}_{F})=\prod_{v\in S_{\infty }}G(F_{v})\times G(\mathbb{A}%
_{f,F})
\end{equation*}%
and so, for any open subset $U$ of $\prod_{v\in S_{\infty }}G(F_{v})$, there
exist elements $\alpha \in G(F_{v})$ and $\beta \in G(F)$ such that $\alpha
\beta \in U\times K$. Clearly $\beta \in \Gamma $, and $\alpha \beta \in U$;
thus $G(F_{v})\Gamma $ is dense in $\prod_{v\in S_{\infty }}G(F_{v})$.
\end{proof}

\begin{remark}
Let $G^{\prime }=\Res_{F/\mathbb{Q}{}}G$ with $G$ and $F$ as in the lemma.
Then the conclusion of the lemma can be restated as follows: let $%
G_{v}^{\prime }$ be a noncompact factor of the Lie group $G^{\prime }(%
\mathbb{R}{})$, and let $\Gamma $ be a congruence subgroup of $G^{\prime }(%
\mathbb{Q}{})$; then $G_{v}^{\prime }\cdot \Gamma $ is dense in $G^{\prime }(%
\mathbb{R}{})$.
\end{remark}

Since $T^{I}(\tilde{X})$ is stable under $\Gamma $, it defines a subbundle $%
T^{I}(X)$ of $T(X)$. By construction, $T^{I}(\tilde{X})$ is involutive, and
so $T^{I}(X)$ is an involutive subbundle of $T(X)$.

\begin{proposition}
Assume $G$ is almost simple over $\mathbb{Q}{}$. If the foliation defined by 
$T^{I}(X)$ on $X$ has a closed leaf, then $I=\emptyset $ or $I=\{1,\ldots
,s\}$.
\end{proposition}

\begin{proof}
Let $Z$ be the closed leaf, and let $\tilde{Z}=p^{-1}(Z)$. Then $\tilde{Z}$
is a closed submanifold of $\tilde{X}$ stable under $\Gamma $ and all $G_{i}(%
\mathbb{R}{})$ for $i\in I$. If $I\neq \emptyset $, then (4.3) shows that $%
\Gamma \cdot \prod G_{i}(\mathbb{R}{})$ is dense in $G(\mathbb{R}{})$, which
acts transitively on $\tilde{X}$. Therefore, $\tilde{Z}=\tilde{X}$ and $Z=X$%
, whence $T^{I}(X)=T(X)$.
\end{proof}

\section{Completion of the proof.}

Let $X$, $\tilde{X}$, $G$, ... be as in \S 1 and assume condition (0.3a),
i.e., that $\codim(\partial X)\geq 3$. Recall the notations 
\begin{equation*}
\begin{diagram} \tilde{X} & \rInto{p} & (\hat{X})^{\text{an}} & & \hat{X} \\
& \rdTo_{p_{K}} & \dTo_{(q_{K})^{\text{an}}} & & \dTo_{q_{K}} \\ & & (X_{K})^{\text{an}} &
& X_{K} \end{diagram} \quad \begin{diagram} \tilde{X}^{\sigma } &
\rInto{p^{\sigma }} & (\sigma \hat{X})^{\text{an}} & & \sigma \hat{X} \\ &
\rdTo_{p_{K}^{\sigma }} & \dTo_{(\sigma q_{K})^{\text{an}}} & & \dTo_{\sigma
q_{K}} \\ & & (\sigma X_{K})^{\text{an}} & & \sigma X_{K} \end{diagram}
\end{equation*}

Let $\mu $ be the Bergmann volume form on $\tilde{X}^{\sigma }$, and let $%
Z_{0}$ be the set on which $\mu $ is zero. As we saw in \S 3, $Z_{0}$ is a
proper subset of $\tilde{X}^{\sigma }$. Clearly, it a complex analytic
subset and is $Q$-invariant. The complement $\tilde{U}_{0}$ of $\tilde{Z}%
_{0} $ in $\tilde{X}^{\sigma }$ is also $Q$-invariant, and (see 3.5) there
is a map $\gamma \colon \tilde{U}_{0}\rightarrow \mathbb{P}{}(\mathcal{H}{}(%
\tilde{X}^{\sigma })^{\vee })$ such that the Bergmann Hermitian form on $%
\tilde{U}_{0}$ is the inverse image of the canonical metric on $\mathbb{P}{}(%
\mathcal{H}{}(\tilde{X}^{\sigma })^{\vee })$. Note that $Q$ acts on $\mathbb{%
P}{}(\mathcal{H}{}(\tilde{X}^{\sigma })^{\vee })$ through its action on $%
\tilde{X}^{\sigma }$, and that $\gamma $ is a $Q$-equivariant map.\footnote{%
Should observe also that $d\gamma $ is not identically zero --- if it were
then $\omega $ would be ``constant'' --- could not have $\int \omega \wedge 
\bar{\omega}<\infty $.} Define 
\begin{equation*}
\tilde{Z}_{1}=\{z\in \tilde{U}_{0}\mid \rank(d\gamma )_{z}<\max_{x}\rank%
(d\gamma )_{x}\}
\end{equation*}

\noindent and $\tilde{Z}=\tilde{Z}_{0}\cup \tilde{Z}_{1}$ --- this is again
a $Q$-invariant complex analytic subset of $\tilde{X}^{\sigma }$ not equal
to $\tilde{X}^{\sigma }$, and so (1.9) shows that $\codim(\tilde{Z})\geq %
\codim(\partial X)\geq 3$. The complement $\tilde{U}=\tilde{X}^{\sigma
}\smallsetminus \tilde{Z}$ of $\tilde{Z}$ is also $Q$-invariant. We have a
diagram 
\begin{diagram}
T(\tilde{U}) & \rTo & T(\mathbb{P}) \\ 
\dTo &  & \dTo \\ 
\tilde{U} & \rTo^{\gamma} & \mathbb{P}
\end{diagram}

\noindent Define 
\begin{equation*}
\tilde{W}=\Ker(T(\tilde{U})\rightarrow \gamma ^{\ast }T(\mathbb{P}{})).
\end{equation*}

\noindent Because it has constant rank, $\tilde{W}$ is a subbundle of $T(%
\tilde{U})$. It is $Q$-invariant.

As in \S 4, let $\tilde{X}=X_{1}\times \cdots \times X_{s}$, and, for each $%
I\subset \{1,\ldots ,s\}$ define $T^{I}(\tilde{X})\subset T(\tilde{X})$ and $%
T^{I}(X)\subset T(X)$. Recall the following result: let $\bar{V}$ be a
complete algebraic variety over $\mathbb{C}{}$, and let $V$ be a nonsingular
subvariety such that $\bar{V}\smallsetminus V$ has codimension $\geq 3$;
then $F\mapsto F^{\text{an}}$ induces an equivalence between the category of
coherent locally free algebraic sheaves on $V$ and that of coherent locally
free analytic sheaves on $V^{\text{an}}$ (see, for example, Hartshorne 1970,
p.~222-223). Thus $T^{I}(X)$ is an algebraic subbundle of $T(X)$, and $%
\sigma T^{I}(X)\subset \sigma T(X)=T(\sigma X)$ is defined. Set $T^{I}(%
\tilde{X}^{\sigma })$ equal to the inverse image of $\sigma T^{I}(X)$ on $%
\tilde{X}^{\sigma }$.

\begin{lemma}
For some $I\subset \{1,2,\ldots ,s\}$, $\tilde{W}=T^{I}(\tilde{X}^{\sigma })|%
\tilde{U}$.
\end{lemma}

\begin{proof}
Since both $\tilde{U}$ and $\tilde{W}$ are $Q$-invariant and $\Gamma
^{\sigma }\subset Q$, we can pass to the quotient and obtain $U=\Gamma
^{\sigma }\backslash \tilde{U}\subset \sigma X$ and $W=\Gamma ^{\sigma
}\backslash \tilde{W}$ a subbundle of $T(U)$. As $\codim(\sigma X^{\ast
}\smallsetminus U)\geq 3$, the result recalled above shows that $W$ is an
algebraic subbundle of $T(U)$. Let $j$ be the inclusion $U\hookrightarrow
\sigma X$. Regard $W$ as a sheaf on $U$, and form $j_{\ast }W$ --- this is a
coherent algebraic sheaf on $\sigma X$. Let $Y$ be the subset of $\sigma X$
where $j_{\ast }W$ is not locally free. Then $Y$ is an algebraic subset of $%
\sigma X$ and its inverse image on $\tilde{X}^{\sigma }$is $Q$-invariant ($%
(p^{\sigma })^{-1}Y$ is the set where $\tilde{j}_{\ast }\tilde{W}$ is not
locally free). Therefore (1.8) shows that $Y$ is empty, and so $j_{\ast }W$
is locally free. A similar argument applied to the support of the kernel of $%
j_{\ast }W\rightarrow T(\sigma X)$ shows that $j_{\ast }W$ is a locally free
subsheaf of $T(\sigma X)$. The fact recalled above shows that it is
algebraic. Similar arguments apply to each variety in the projective system $%
\sigma \hat{X}$, and so we obtain an algebraic subbundle $\hat{W}\subset
T(\sigma \hat{X})$ invariant under $Q$, and therefore under $G(\mathbb{A}%
_{f})$. Hence $\sigma ^{-1}\hat{W}$ is a subbundle of $T(\hat{X})$ invariant
under $G(\mathbb{A}_{f})$, and $(\sigma ^{-1}\hat{W})^{\text{an}}|\tilde{X}$
is $G(\mathbb{R}{})$-invariant. Now (4.1) shows that $(\sigma ^{-1}\hat{W})^{%
\text{an}}|\tilde{X}=T^{I}(\tilde{X})$ for some $I$ and, because $G(\mathbb{A%
}_{f})$ acts transitively on the set of connected components of $\hat{X}^{%
\text{an}}$, this implies that $\sigma ^{-1}\hat{W}=T^{I}(\hat{X})$.
\end{proof}

The condition that a subbundle of the tangent bundle to an algebraic variety
be involutive is algebraic. Thus $T^{I}(\tilde{X}^{\sigma })$ is involutive.

\begin{lemma}
The foliation of $U=\Gamma ^{\sigma }\backslash \tilde{U}$ defined by $%
W=\Gamma ^{\sigma }\backslash \tilde{W}$ has a closed leaf.
\end{lemma}

\begin{proof}
The reader is invited to check for himself the proof of Kazhdan (Kazhdan 1983,
p.~153--156) --- essentially the leaves of the foliation are the equivalence
classes for the relation defined in p.~153. Alternatively, it may be
possible to give a proof along the following lines: as we observed above,
there is a $Q$-equivariant map $\gamma \colon \tilde{U}\rightarrow \mathbb{P}%
{}(\mathcal{H}{}(\tilde{X}^{\sigma })^{\vee })$; if we can pass to the
quotient by $\Gamma ^{\sigma }$, we get a map $U=\Gamma ^{\sigma }\backslash 
\tilde{U}\rightarrow \Gamma ^{\sigma }\backslash \mathbb{P}{}(\mathcal{H}{}(%
\tilde{X}^{\sigma })^{\vee })$ whose fibres are the leaves of the foliation.
\end{proof}

We now assume also the conditions (0.3b). Thus $G=\Res_{F/\mathbb{Q}%
{}}G^{\prime }$ with $G^{\prime }$ absolutely simple and $F$ totally real,
and for some special point $\tilde{x}\in \tilde{X}$, the maximal torus $%
T\subset G$ fixing $\tilde{x}$ is of the form $\Res_{F/\mathbb{Q}%
{}}T^{\prime }$ where $T^{\prime }$ splits over a quadratic imaginary
extension $L$ of $F$.

As $T_{L}^{\prime }$ is split, we can write 
\begin{equation*}
\mathfrak{g}{}_{L}^{\prime }=\mathfrak{t}{}_{L}^{\prime }\oplus
\bigoplus_{\alpha \in R}(\mathfrak{g}{}_{L}^{\prime })_{\alpha },\quad 
\mathfrak{g}{}_{L}^{\prime }=\Lie(G_{L}^{\prime }),\quad \mathfrak{t}%
{}_{L}^{\prime }=\Lie(T_{L}^{\prime }),
\end{equation*}

\noindent where $R$ is the set of roots $R(G_{\mathbb{\mathbb{C}{}}}^{\prime
},T_{\mathbb{C}{}}^{\prime })$. For each $\alpha \in R$, let 
\begin{equation*}
\mathfrak{h}{}_{\alpha }^{\prime }=\mathfrak{t}{}_{L}^{\prime }\oplus (%
\mathfrak{g}{}_{L}^{\prime })_{\alpha }\oplus (\mathfrak{g}{}_{L}^{\prime
})_{-\alpha }.
\end{equation*}

\noindent It is defined over $F$, and we let $H_{\alpha }^{\prime }$ be the
corresponding connected subgroup of $G^{\prime }$. Let $H_{\alpha }=\Res_{F/%
\mathbb{Q}{}}H_{\alpha }^{\prime }$ --- it is a reductive group of type $%
A_{1}$ and $H_{\alpha }^{\text{ad}}$ is $\mathbb{Q}{}$-simple.

Let $\tilde{X}_{\alpha }$ be the orbit of $\tilde{x}$ under the action of $%
H_{\alpha }(\mathbb{R}{})$. Then $\tilde{X}_{\alpha }$ is a Hermitian
symmetric domain (possibly consisting of one element) with 
\begin{equation*}
\Aut(\tilde{X}_{\alpha })^{+}=(H_{\alpha }(\mathbb{R}{})/\{\text{maximal
compact normal subgroup}\})^{+}\text{.}
\end{equation*}

\noindent It follows from Deligne 1971, 1.15, that the projective system of
arithmetic varieties $\hat{X}_{\alpha }$ embeds into $\hat{X}$. We assume
Theorem 0.1 for the families $\hat{X}_{\alpha }$; in particular, $\hat{X}%
_{\alpha }$ is the family associated with a $\mathbb{Q}{}$-group $H_{\alpha
}^{\sigma }$.

\begin{lemma}
The tangent space $T_{\tilde{x}}(\tilde{X})$ is generated by the subspaces $%
T_{\tilde{x}}(\tilde{X}_{\alpha })$. (Note, $\tilde{x}$ is as defined above.)
\end{lemma}

\begin{proof}
Easy.
\end{proof}

\begin{lemma}
Let $\tilde{x}^{\sigma }$be any point of $\tilde{X}^{\sigma }$such that $%
\tilde{x}^{\sigma }=g(\sigma \hat{x})$ for some $g\in G(\mathbb{A}_{f})$,
where $\hat{x}$ is the image of $\tilde{x}$ (see above) in $\hat{X}$. Then $%
\tilde{x}^{\sigma }\in \tilde{U}$.
\end{lemma}

\begin{proof}
Otherwise $\tilde{Z}\supset Q\tilde{x}^{\sigma }\supset H_{\alpha }^{\sigma
}(\mathbb{Q}{})\tilde{x}^{\sigma }$, and so $\tilde{Z}\supset \tilde{X}%
_{\alpha }^{\sigma }$ for all $\alpha $. This implies that $\dim \tilde{Z}%
=\dim \tilde{X}^{\sigma }$, which contradicts an earlier statement.
\end{proof}

Let 
\begin{equation*}
J=\{\tau \colon F\rightarrow \mathbb{R}{}\mid G_{\tau }^{\prime }\text{ is
not compact}\}\text{,}
\end{equation*}
where $G_{\tau }^{\prime }=G\otimes _{F,\tau }\mathbb{R}{}$. Then $J$
indexes the irreducible components of $\tilde{X}$, and we shall use it for
this purpose rather than $\{1,\ldots ,s\}$. Thus $I$ is now a subset of $J$.
Let 
\begin{equation*}
J_{\alpha }=\{\tau \colon F\rightarrow \mathbb{R}{}\mid H_{\alpha ,\tau
}^{\prime }\text{ is not compact}\}\text{.}
\end{equation*}

\begin{lemma}
Either $J_{\alpha }\cap I=\emptyset $ or $J_{\alpha }\subset I$.
\end{lemma}

\begin{proof}
The set $\tilde{Z}\cap \tilde{X}_{\alpha }^{\sigma }$ is stable under $Q\cap
H_{\alpha }(\mathbb{A}_{f})=H_{\alpha }(\mathbb{Q}{})/Z_{\alpha }(\mathbb{Q}%
{})$ and so is either empty of all of $\tilde{X}_{\alpha }^{\sigma }$. The
last lemma shows that it is not equal to $\tilde{X}_{\alpha }^{\sigma }$,
and so we have that $\tilde{X}_{\alpha }^{\sigma }\subset \tilde{U}$.

From (4.1), (5.2), and (4.4), we know that $W|\sigma X_{\alpha }$ is $0$ or $%
T(\sigma X_{\alpha })$ (recall that we are assuming $\sigma X_{\alpha }$ is
arithmetic). Thus $\sigma ^{-1}W$ on $X_{\alpha }$ is $0$ or $T(X_{\alpha })$%
. But (by 5.1), $\sigma ^{-1}W=T^{I}(X)$ for some $I\subset J$. We conclude
that $I\supset J_{\alpha }$ (and $\sigma ^{-1}W|X_{\alpha }=T(X_{\alpha })$)
or $I\cap J_{\alpha }=\emptyset $ (and $\sigma ^{-1}W|X_{\alpha }=0$).
\end{proof}

Note that it is not possible for $I$ to contain all $J_{\alpha }$, for then $%
I=J$ and the Bergmann Hermitian form is identically zero (cf.~the first page of
this section). On the other hand, if $I\cap J_{\alpha }=\emptyset $, then $%
I= $ $\emptyset $, $\tilde{W}=0$, and we can apply (2.13) to complete the
proof of (0.1).

Thus it remains to show that the hypothesis that $I\cap J_{\alpha
}=\emptyset $ for some, but not all, $\alpha $ leads to a contradiction.

\begin{lemma}
Let $\mathfrak{h}{}_{I}^{\prime }=\sum_{J_{\alpha }\subset I}\mathfrak{h}%
{}_{\alpha }^{\prime }$; then $\mathfrak{h}{}_{I}^{\prime }$ is a
sub-Lie-algebra of $\mathfrak{g}{}^{\prime }$.
\end{lemma}

\begin{proof}
Let $\mathfrak{g}{}_{\tau }^{\prime }=\mathfrak{k}{}_{\tau }^{\prime }\oplus 
\mathfrak{p}{}_{\tau }^{\prime }$ as usual. Then 
\begin{equation*}
J_{\alpha }\subset I\iff \text{for all }\tau \notin I\text{, }\mathfrak{h}%
{}_{\alpha ,\tau }^{\prime }\subset \mathfrak{k}{}_{\tau }^{\prime }.
\end{equation*}

As $\mathfrak{k}{}_{\tau }^{\prime }$ is a subalgebra of $\mathfrak{g}%
{}_{\tau }^{\prime }$, the result is now obvious.
\end{proof}

Let $H_{I}^{\prime }$ be the subgroup of $G^{\prime }$ corresponding to $%
\mathfrak{h}{}_{I}^{\prime }$. Then:

\begin{itemize}
\item for $\tau \notin I$, $H_{I,\tau }^{\prime }$ is compact (and hence $%
H_{I}^{\prime }$ is reductive),

\item for $\tau \in I$, $H_{I,\tau }^{\prime }(\mathbb{R}{})\tilde{x}_{\tau
}\supset (\tilde{X}_{\alpha })_{\tau }$ for all $\alpha $, and so $H_{I,\tau
}^{\prime }(\mathbb{R}{})\tilde{x}_{\tau }=\tilde{X}_{\tau }$ (apply
Kobayashi and Nomizu, 1963, p.~178, 4.8 if it isn't obvious).
\end{itemize}

\noindent (We are writing $\tilde{x}=(\ldots ,\tilde{x}_{\tau },\ldots )\in
\cdots \times \tilde{X}_{\tau }\times \cdots $.) This shows that $H_{I,\tau
}^{\prime }=G_{\tau }^{\prime }$ (see 2.3), and so $H_{I}^{\prime
}=G^{\prime }$. This is the contradiction that proves the theorem.

\section*{References}

\setlength{\parindent}{-1em}

{\footnotesize Borel, Armand, K\"ahlerian coset spaces of semisimple Lie
groups. Proc. Nat. Acad. Sci. U. S. A. 40 (1954), 1147--1151. }

{\footnotesize Bourbaki, N., General Topology, Chapters 1--4, Springer,
1989. }

{\footnotesize Deligne, Pierre, Travaux de Shimura. S\'eminaire Bourbaki,
23\`eme ann\'{e}e (1970/71), Exp. No. 389, pp. 123--165. Lecture Notes in
Math., Vol. 244, Springer, Berlin, 1971. }

{\footnotesize Deligne, Pierre, Vari\'et\'es de Shimura: interpr\'etation
modulaire, et techniques de construction de mod\`eles canoniques.
Automorphic forms, representations and $L$-functions (Proc. Sympos. Pure
Math., Oregon State Univ., Corvallis, Ore., 1977), Part 2, pp. 247--289,
Proc. Sympos. Pure Math., XXXIII, Amer. Math. Soc., Providence, R.I., 1979. }

{\footnotesize Deligne, Pierre, Hodge cycles on abelian varieties (Notes by
J.S. Milne). In Deligne, Pierre; Milne, James S.; Ogus, Arthur; Shih,
Kuang-yen, Hodge cycles, motives, and Shimura varieties. Lecture Notes in
Mathematics, 900. Springer-Verlag, Berlin-New York, 1982, pp.\ 9--100. }

{\footnotesize Deligne, Pierre, Motifs et groupe de Taniyama, in Deligne,
Pierre; Milne, James S.; Ogus, Arthur; Shih, Kuang-yen, Hodge cycles,
motives, and Shimura varieties. Lecture Notes in Mathematics, 900.
Springer-Verlag, Berlin-New York, 1982, pp.\ 261--279. }

{\footnotesize Griffiths, Phillip A., Entire holomorphic mappings in one and
several complex variables. The fifth set of Hermann Weyl Lectures, given at
the Institute for Advanced Study, Princeton, N. J., October and November
1974. Annals of Mathematics Studies, No. 85. Princeton University Press,
Princeton, N. J.; University of Tokyo Press, Tokyo, 1976. }

{\footnotesize Hano, Jun-ichi, On Kaehlerian homogeneous spaces of
unimodular Lie groups. Amer. J. Math. 79 1957 885--900. }

{\footnotesize Hartshorne, Robin, Ample subvarieties of algebraic varieties.
Notes written in collaboration with C. Musili. Lecture Notes in Mathematics,
Vol. 156 Springer-Verlag, Berlin-New York 1970. }

{\footnotesize Helgason, S., Differential geometry and symmetric spaces.
Pure and Applied Mathematics, Vol. XII. Academic Press, New York-London
1962. }

{\footnotesize Hironaka, Heisuke, Resolution of singularities of an
algebraic variety over a field of characteristic zero. I, II. Ann. of Math.
(2) 79 (1964), 109--203; ibid. (2) 79 (1964), 205--326. }

{\footnotesize Kazhdan, David, Arithmetic varieties and their fields of
quasi-definition. Actes du Congr\`{e}s International des Math\'{e}maticiens
(Nice, 1970), Tome 2, pp. 321--325. Gauthier-Villars, Paris, 1971. }

{\footnotesize Kazhdan, David, On arithmetic varieties. Lie groups and their
representations (Proc. Summer School, Bolyai J\'{a}nos Math. Soc., Budapest,
1971), pp. 151--217. Halsted, New York, 1975. }

{\footnotesize Kazhdan, David, On arithmetic varieties. II. Israel J. Math.
44 (1983), no. 2, 139--159. }

{\footnotesize Kobayashi, Shoshichi, Geometry of bounded domains. Trans.
Amer. Math. Soc. 92 (1959) 267--290. }

{\footnotesize Kobayashi, Shoshichi; Nomizu, Katsumi, Foundations of
differential geometry. Vol I. Interscience Publishers, a division of John
Wiley, New York-London 1963. }

{\footnotesize Kobayashi, Shoshichi; Nomizu, Katsumi, Foundations of
differential geometry. Vol. II. Interscience Tracts in Pure and Applied
Mathematics, No. 15 Vol. II Interscience Publishers John Wiley Inc., New
York-London-Sydney 1969. }

{\footnotesize Kobayashi, Shoshichi, Transformation groups in differential
geometry. Ergebnisse der Mathematik und ihrer Grenzgebiete, Band 70.
Springer-Verlag, New York-Heidelberg, 1972. }

{\footnotesize Koszul, J. L., Sur la forme hermitienne canonique des espaces
homog\`{e}nes complexes. Canad. J. Math. 7 (1955), 562--576. }

{\footnotesize Koszul, J.-L., Expos\'{e}s sur les espaces homog\`{e}nes sym%
\'{e}triques. Publica\c{c}ao da Sociedade de Matematica de Sao Paulo, Sao
Paulo 1959 71 pp. }

{\footnotesize Langlands, R. P., Automorphic representations, Shimura
varieties, and motives. Ein M\"{a}rchen. Automorphic forms, representations
and $L$-functions (Proc. Sympos. Pure Math., Oregon State Univ., Corvallis,
Ore., 1977), Part 2, pp. 205--246, Proc. Sympos. Pure Math., XXXIII, Amer.
Math. Soc., Providence, R.I., 1979. }

{\footnotesize Margulis, G.A., Arithmeticity of irreducible lattices in
semisimple groups of rank greater than 1, AMS Transl. 109 (1977), 33-45. }

{\footnotesize Milne, J.S., and Shih, K-y., Conjugates of Shimura varieties,
Deligne, Pierre; Milne, James S.; Ogus, Arthur; Shih, Kuang-yen Hodge
cycles, motives, and Shimura varieties. Lecture Notes in Mathematics, 900.
Springer-Verlag, Berlin-New York, 1982, pp.~280--356.}

{\footnotesize Milne, J. S., The action of an automorphism of $\mathbf{C}$
on a Shimura variety and its special points. Arithmetic and geometry, Vol.
I, 239--265, Progr. Math., 35, Birkh\"{a}user Boston, Boston, MA, 1983. }

{\footnotesize Narasimhan, Raghavan, Introduction to the theory of analytic
spaces. Lecture Notes in Mathematics, No. 25 Springer-Verlag, Berlin-New
York 1966. }

{\footnotesize Piateski-Shapiro, I. I., Automorphic functions and the
geometry of classical domains. Translated from the Russian. Mathematics and
Its Applications, Vol. 8 Gordon and Breach Science Publishers, New
York-London-Paris 1969 }

{\footnotesize Raghunathan, M. S., Discrete subgroups of Lie groups.
Ergebnisse der Mathematik und ihrer Grenzgebiete, Band 68. Springer-Verlag,
New York-Heidelberg, 1972. }

{\footnotesize Yau, Shing Tung, On the Ricci curvature of a compact K\"ahler
manifold and the complex Monge-Amp\`ere equation. I. Comm. Pure Appl. Math.
31 (1978), no. 3, 339--411. }

{\normalsize\vfill\hfill \TeX ed July 12, 2001. }

\end{document}